\documentclass[11pt]{article}
\usepackage{amssymb,amsmath,amsthm}
\usepackage{amscd}
\usepackage{mathrsfs}
\allowdisplaybreaks[1]
\newcommand{\C}{\mathbb{C}}
\newcommand{\N}{\mathbb{N}}
\newcommand{\R}{\mathbb{R}}

\newcommand{\dd}{\,{\rm d}}

\numberwithin{equation}{section}
\newtheorem{thm}{Theorem}[section]
\newtheorem{df}[thm]{Definition}
\newtheorem{prop}[thm]{Proposition}
\newtheorem{lem}[thm]{Lemma}
\newtheorem{rem}[thm]{Remark}
\newtheorem{cor}[thm]{Corollary}

\begin{document}

\title{On Poisson operators and Dirichlet-Neumann maps in $H^s$ for divergence form elliptic operators with Lipschitz coefficients}

\date{}

\author{
\null\\
Yasunori Maekawa\\    
Mathematical Institute, Tohoku University\\          
6-3 Aoba, Aramaki, Aoba, Sendai 980-8578, Japan\\
{\tt maekawa@math.tohoku.ac.jp }
\and
\\
Hideyuki Miura \\
Department of Mathematics,  Graduate School of Science,  Osaka University\\
1-1 Machikaneyama, Toyonaka, Osaka 560-0043, Japan\\
{\tt miura@math.sci.osaka-u.ac.jp}}

\maketitle

\begin{center}
{\bf Abstract}
\end{center}
\vspace{-2mm}

We consider second order uniformly elliptic operators of divergence form in $\R^{d+1}$ whose coefficients are independent of one variable. 
Under the Lipschitz condition on the coefficients 
we characterize the domain of 
the Poisson operators and the Dirichlet-Neumann maps 
in the Sobolev space $H^s(\R^d)$ for each $s\in [0,1]$. Moreover, 
we also show a factorization formula for the elliptic operator
in terms of the Poisson operator.

\vspace{0.5cm}


               
\section{Introduction}\label{sec.intro}


In this present paper we consider the second order elliptic operator of divergence form in  $\R^{d+1}= \{ (x,t) \in \R^d \times \R \}$, 
\begin{equation}
\mathcal{A} = -\nabla \cdot A \nabla, ~~~~~~~~~~~~~~A = A (x) = \big (a_{i,j} (x) \big ) _{1\leq i,j\leq d+1}~.\label{def.calA}
\end{equation}
Here $d\in \N$, $\nabla = (\nabla_x,\partial_t)^\top$ with $\nabla_x=(\partial_1,\cdots,\partial_d)^\top$,  and each $a_{i,j}$ is complex-valued and  assumed to be $t$-independent. The adjoint matrix of $A$ will be denoted by $A^*$. We assume the uniformly ellipticity condition
\begin{align}
{\rm Re} \langle A (x) \eta, \eta \rangle \geq \nu_1 |\eta|^2, ~~~~~~~~~~ | \langle A (x) \eta,\zeta\rangle | \leq \nu_2 |\eta| |\zeta|\label{ellipticity}
\end{align}
for all $\eta,\zeta\in \C^{d+1}$ with positive constants  $\nu_1,\nu_2$. 
Here $\langle \cdot,\cdot \rangle$ denotes the inner product of $\C^{d+1}$, i.e., $\langle \eta, \zeta\rangle = \sum_{j=1}^{d+1}\eta_j \bar{\zeta}_j$ for $\eta,\zeta\in \C^{d+1}$. For later use we set 
\begin{align*}
A' =(a_{i,j})_{1\leq i,j\leq d},~~~b=a_{d+1, d+1}, ~~~{\bf r_1} = ( a_{1, d+1},\cdots , a_{d, d+1} )^\top, ~~~{\bf r_2}  = ( a_{d+1,1} ,\cdots , a_{d+1,d})^\top.
\end{align*}
We will also use the notation $\mathcal{A}'=-\nabla_x\cdot A' \nabla_x$. 
In this paper, we are concerned with the Poisson operator and 
the Dirichlet-Neumann map associated with $\mathcal{A}$, 
which play fundamental roles in the boundary value problems
for the elliptic operators.
They are defined through $\mathcal{A}$-extension of the boundary data on $\R^d=\partial\R^{d+1}_+$ to the upper half space.
\begin{df}\label{def.A-extension} {\rm (i)} For a given $h\in \mathcal{S}' (\R^d)$ we denote by $M_h: \mathcal{S}(\R^d)\rightarrow \mathcal{S}'(\R^d)$ the multiplication $M_h u = h u$.
 
\noindent {\rm (ii)} We denote by $E_{\mathcal{A}}: \dot{H}^{1/2} (\R^d)\rightarrow \dot{H}^1 (\R^{d+1}_+)$ the $\mathcal{A}$-extension operator, i.e., $w=E_{\mathcal{A}} f$ is the solution to the Dirichlet problem
\begin{equation}\label{eq.dirichlet0}
\begin{cases}
& \mathcal{A} u  = 0~~~~~~~{\rm in}~~~\R^{d+1}_+,\\
& \hspace{0.3cm} u  = f~~~~~~~{\rm on}~~\partial\R^{d+1}_+=\R^d.
\end{cases}
\end{equation}
The one parameter family of linear operators $\{E_{\mathcal{A}} (t)\}_{t\geq 0}$, defined by $E_{\mathcal{A}} (t) f = (E_{\mathcal{A}} f )(\cdot,t)$ for $f\in \dot{H}^{1/2}(\R^d)$, is called the Poisson semigroup associated with $\mathcal{A}$.

\noindent {\rm (iii)} We denote by $\Lambda_{\mathcal{A}}: D_{L^2} (\Lambda_{\mathcal{A}}) \subset \dot{H}^{1/2} (\R^d)\rightarrow \dot{H}^{-1/2} (\R^d)$ the Dirichlet-Neumann map associated with $\mathcal{A}$, which is defined through the sesquilinear form 
\begin{equation}
\langle  \Lambda_{\mathcal{A}} f, g\rangle_{\dot{H}^{-\frac12},\dot{H}^{\frac12}}  =  \langle  A\nabla E_{\mathcal{A}}f, \nabla E_{\mathcal{A}} g \rangle_{L^2(\R^{d+1}_+)},~~~~~~~~~~f,g \in \dot{H}^{\frac12}(\R^d).\label{def.Lambda}
\end{equation}
Here $\langle \cdot,\cdot\rangle _{\dot{H}^{-1/2},\dot{H}^{1/2}}$ denotes the duality coupling of  $\dot{H}^{-1/2}(\R^d)$ and $\dot{H}^{1/2}(\R^d)$. 

\end{df}
\noindent Here $\dot{H}^s (\R^d)$ is the homogeneous Sobolev space of the order $s\in \R$ and $D_H(T)$ denotes the domain of a linear operator $T$ 
in a Banach space $H$. Since the ellipticity condition \eqref{ellipticity}
ensures that $E_\mathcal{A}$ is well-defined in $\dot{H}^{1/2}(\R^d)$
via the Lax-Milgram theorem, 
it is not difficult to see that $\{E_{\mathcal{A}} (t)\}_{t\geq 0}$ 
is realized as a strongly continuous and analytic semigroup 
in $\dot{H}^{1/2}(\R^d)$ and in $H^{1/2}(\R^d)$ (see, e.g. \cite[Proposition 2.4]{MaekawaMiura1}). 
Then the generator of the Poisson semigroup will be denoted by $-\mathcal{P}_{\mathcal{A}}$, and $\mathcal{P}_{\mathcal{A}}$ is called the {\it Poisson operator} (associated with $\mathcal{A}$). 
As for the Dirichlet-Neumann map, it is well known from the theory of sesquilinear forms that  \eqref{ellipticity} guarantees 
the generation of a strongly continuous and analytic semigroup in $L^2 (\R^d)$; see \cite{Kato}. On the other hand, the realization of the Poisson semigroup in $L^2 (\R^d)$ is nothing but the solvability of the elliptic boundary 
value problem \eqref{eq.dirichlet0} for $L^2$ boundary data (see \cite{MaekawaMiura1} for details), there have been a lot of works on 
this subject by now.  Moreover the characterization of $D_{L^2}(\mathcal{P}_{\mathcal{A}})$ is studied as well, for it provides precise informations on 
the behavior of $\mathcal{A}$-extension near the boundary.
As far as the authors know, 
these problems are affirmatively settled at least for the following classes of $A$. 

\vspace{0.2cm}

 (I) $A$ is a constant matrix, i.e., $A(x)=A$;  (II) $A$ is Hermite, i.e., $A^*=A$; (III) $A$ is block type, i.e., ${\bf r_1}={\bf r_2} =0$; (IV) $A$ is a small $L^\infty$ perturbation of $B$ satisfying one of (I)-(III) above.

\vspace{0.2cm}

The case (I) is easy since one can directly derive the solution formula for  
\eqref{eq.dirichlet0} with the aid of the Fourier transform. The case (II) is a classical problem, for it is closely related with the Laplace equations in  Lipschitz domains,
and it is studied in \cite{Dahlberg1, JerisonKenig1, JerisonKenig2,
 Verchota, Dahlberg2, KenigPipher, Aucher.et.al.2, Aucher.et.al.3}. 
The case (III) is considered in 
\cite{Aucher.et.al.1, Aucher.et.al.2}. In this case, the Poisson operator 
essentially coincides with the Dirichlet-Neumann map, and the characterization 
$D_{L^2}(\mathcal{P}_{\mathcal{A}})=H^1 (\R^d)$ 
is known as the Kato square root problem 
for divergence form elliptic operators, 
which is settled by \cite{Aucher.et.al.1}. 
The case (IV) is solved in \cite{Fabes.et.al.1} 
when $B$ is a constant matrix, 
and in \cite{Aucher.et.al.2,Aucher.et.al.3,Alfonseca.et.al.1} when $B$ is a Hermite, or block matrix. 
Recently in \cite{MaekawaMiura1}, the authors of the present
paper showed the $L^2$ solvablity of 
\eqref{eq.dirichlet0} and verified 
the characterization $D_{L^2}(\mathcal{P}_{\mathcal{A}})=H^1 (\R^d)$, 
when ${\bf r_1}$, ${\bf r_2}$, and $b$ are real, and 
 $\nabla_x \cdot({\bf r_1}+{\bf r_2})$ 
belong to $L^d(\R^d)+L^\infty(\R^d)$.
We note that in the cases (II)-(IV) 
the coefficients of $A$ are not discontinuous in general.
However, it is shown in \cite{Kenig.et.al.1} 
that if one imposes only \eqref{ellipticity} and the coefficients are 
discontinuous, 
the Dirichlet problem \eqref{eq.dirichlet0} is not always 
solvable for boundary data in $L^2 (\R^d)$.
This means that some additional conditions on $A$ such as 
(I)-(IV) are required 
in order to extend the Poisson semigroup in $H^{1/2}(\R^d)$ 
as a semigroup in $L^2 (\R^d)$. 

As our first result, 
we show the realization of the Poisson semigroup 
and the characterization of the domain of the generator 
 in $H^s(\R^d)$ for $s \in [0,1]$ under the Lipschitz regularity assumption:
\begin{align}
{\rm Lip} (A) = \sum_{i,j} \sup_{x,y\in \R^d} \frac{|a_{i,j} (x) - a_{i,j} (y) |}{|x-y|} <\infty. \label{lipschitz}
\end{align}
The precise statement of the result is given as follows:  
\begin{thm}\label{thm.main.1} Let $A=A(x)$ be a $t$-independent complex coefficient matrix satisfying \eqref{ellipticity} and \eqref{lipschitz}. Then the following statements hold:
\\
{\rm (i)} The Poisson semigroup  in $H^{1/2}(\R^d)$ is extended as a strongly continuous and analytic semigroup in $L^2 (\R^d)$. Moreover, its generator $-\mathcal{P}_{\mathcal{A}}$  satisfies $D_{L^2}(\mathcal{P}_{\mathcal{A}}) = H^1 (\R^d)$ with equivalent norms, and $\mathcal{P}_{\mathcal{A}}$ admits a bounded $H^\infty$ calculus  in $L^2 (\R^d)$.

\noindent {\rm (ii)} Let $s\in [0,1]$. Then $H^s(\R^d)$ is invariant under the action of the Poisson semigroup $\{e^{-t\mathcal{P}_{\mathcal{A}}}\}_{t\geq 0}$ in $L^2 (\R^d)$, and its restriction on $H^s(\R^d)$ defines a strongly continuous and analytic semigroup in $H^s(\R^d)$. Moreover, its generator, denoted again by $-\mathcal{P}_{\mathcal{A}}$, satisfies $D_{H^s}(\mathcal{P}_{\mathcal{A}}) = H^{1+s}(\R^d)$ with equivalent norms.  

\end{thm}

For the definition of bounded $H^\infty$ calculus for sectorial operators, see, e.g., \cite[Chapter 5]{Haase}. 
The main feature in this result is that we do not assume any {\it structural} conditions such as (I)-(IV). 
In contrast to approaches taken in the aforementioned results, 
we analyze $\mathcal{P}_{\mathcal{A}}$ by looking 
at its principal symbol, which is explicitly calculated as 
\begin{align}
\mu_{\mathcal{A}} (x,\xi) = - \frac{{\bf v} (x) \cdot \xi}{2} + i \big \{ \frac{1}{b (x) } \langle A'(x) \xi,\xi\rangle -\frac{1}{4}({\bf v}(x) \cdot \xi )^2\big \}^\frac12,~~~~~~~~{\bf v} = \frac{{\bf r_1} + {\bf r_2}}{b}.\label{mu_A.intro}
\end{align}
Here  $x,\xi\in \R^d$. 
As is expected, the associated pseudo-differential operator 
$-i\mu_{\mathcal{A}}(\cdot,D_x)$ is shown to be an 
approximation of $\mathcal{P}_{\mathcal{A}}$. 
Since $\mu_{\mathcal{A}}$ is Lipschitz in $x$ 
and homogeneous of degree $1$ in $\xi$, 
one may apply the general theory of pseudo-differential operators with nonsmooth symbols to 
$\mu_{\mathcal{A}}(\cdot,D_x)$; \cite{KumanogoNagase, Marschall, Taylor, Abels, ES} to show Theorem \ref{thm.main.1} at least for $s <1$. 
Here we provide another approach to 
the analysis of $\mathcal{P}_{\mathcal{A}}$ 
which does not rely on the detailed properties of 
$\mu_{\mathcal{A}}(\cdot,D_x)$ obtained from this general theory; 
see Remark \ref{rem.pseudo}. 
Indeed, the key ingredient underlying the proof of Theorem \ref{thm.main.1} 
in our argument is the factorizations 
of operators $\mathcal{A}'$ and $\mathcal{A}$ 
in terms of $\mathcal{P}_{\mathcal{A}}$, which 
we will state in the next theorem. 
We note here that the assertion (ii) includes the critical case $s=1$, 
which seems to be out of reach of general theory of the pseudo-differential 
operators in the works cited above.

\begin{rem}
{\rm
If $\R^{d+1}_+$ is replaced by a bounded Lipschitz domain 
satisfying $VMO$ conditions for the unit normal of the boundary, 
 then the Dirichlet problem  for $A$ with $VMO$ coefficients is 
solved by \cite{MMS} in $L^p$ and Besov spaces. 
 In view of local regularity, 
the Lipschitz condition \eqref{lipschitz} 
assumed in our paper is rather strong. 
However, one has to be careful about the lack of the compactness
of the boundary in our case. In fact, the authors in \cite{MMS} apply 
 a localization argument which enables them to 
approximate $A(x)$ by a constant matrix in each localized domain, then 
 the boundary value problem is reduced to a finite sum of the problems 
for small $VMO$ perturbation of the constant matrices.
However, in our case one cannot use such a localization procedure, since
$A(x)$ is not necessarily close to a constant matrix 
as $|x|\rightarrow \infty$. This difficulty is overcome 
with the aid of the calculus of the symbol \eqref{mu_A.intro}. 
}
\end{rem}

In order to state the next result, 
let us recall the realization of 
$\mathcal{A}$ in $L^2 (\R^{d+1})$:
\begin{align}
D_{L^2}(\mathcal{A}) & = \big \{ u\in H^1 (\R^{d+1})~|~{\rm there ~is}~ F\in L^2 (\R^{d+1})~{\rm such ~that}~\nonumber \\
& ~~~~~~~~~~~~~~~~~~ \langle A\nabla u, \nabla v\rangle _{L^2(\R^{d+1})} = \langle F, v\rangle _{L^2 (\R^{d+1})}~{\rm for~all}~v\in H^1 (\R^{d+1})\big \},\label{realization.A}\\
 \mathcal{A} u & = F ~~~~~ {\rm for} ~~ u\in D_{L^2}(\mathcal{A}).\nonumber 
\end{align}
Note that $D_{L^2}(\mathcal{A}) = H^2 (\R^{d+1})$ holds with equivalent norms because of \eqref{lipschitz}. The realization of $\mathcal{A}' = - \nabla_x\cdot A'\nabla_x$ in $L^2 (\R^d)$ is defined in the similar manner, and we have $D_{L^2}(\mathcal{A}') = H^2 (\R^d)$ with equivalent norms since $A'$ is Lipschitz continuous. 
The following theorem shows the factorization of operators $\mathcal{A}$ 
and $\mathcal{A'}$, and it clarifies the relation of the Poisson operator and the Dirichlet-Neumann map: 
\begin{thm}\label{thm.main.2} Under the same assumption as in Theorem 
\ref{thm.main.1}, the following statements hold:
\\
{\rm (i)} The realization of $\mathcal{A}'$ in $L^2 (\R^d)$ and the realization of $\mathcal{A}$ in $L^2 (\R^{d+1})$ are respectively factorized as
\begin{align}
\mathcal{A}' & = M_b \mathcal{Q}_{\mathcal{A}} \mathcal{P}_{\mathcal{A}}, ~~~~~~~~~\mathcal{Q}_{\mathcal{A}} = M_{1/b} ( M_{\bar{b}} \mathcal{P}_{\mathcal{A}^*} )^*,\label{eq.thm.main.2.1}\\
\mathcal{A}  & = - M_b (\partial_t - \mathcal{Q}_{\mathcal{A}}) (\partial_t + \mathcal{P}_{\mathcal{A}}).\label{eq.thm.main.2.2}
\end{align}
Here  $( M_{\bar{b}} \mathcal{P}_{\mathcal{A}^*} )^*$ is the adjoint operator of $M_{\bar{b}} \mathcal{P}_{\mathcal{A}^*}$ in $L^2 (\R^d)$, while $\mathcal{P}_{\mathcal{A}^*}$ is the Poisson operator in $L^2 (\R^d)$ associated with $\mathcal{A}^*=-\nabla\cdot A^* \nabla$.

\noindent {\rm (ii)} It follows that $D_{L^2}(\Lambda_{\mathcal{A}}) = D_{L^2}(\mathcal{Q}_{\mathcal{A}}) = H^1 (\R^d)$ with equivalent norms and that 
\begin{align}
\mathcal{P}_{\mathcal{A}} & = M_{1/b} \Lambda_{\mathcal{A}} + M_{{\bf r_2}/b} \cdot \nabla_x, \label{eq.thm.main.2.3} \\
\mathcal{Q}_{\mathcal{A}} & = M_{1/b} \Lambda_{\mathcal{A}}  - M_{{\bf r_1}/b}  \cdot \nabla_x  -  M_{(\nabla_x\cdot {\bf r_1})/b}, \label{eq.thm.main.2.4}
\end{align}
as the operators in $L^2 (\R^d)$. 

\end{thm}

Finally we state the counterpart of Theorem \ref{thm.main.1} for
the Dirichlet-Neumann map in $H^s(\R^d)$.

\begin{thm}\label{thm.main.3} 
Let $s \in [0,1]$. Under the same assumption as in Theorem \ref{thm.main.1}, $H^s(\R^d)$ is invariant under the action of  the Dirichlet-Neumann semigroup $\{e^{-t\Lambda_{\mathcal{A}}}\}_{t\geq 0}$ in $L^2 (\R^d)$, and its restriction on $H^s(\R^d)$ defines a strongly continuous and analytic semigroup in $H^s(\R^d)$. Moreover, its generator, denoted again by $-\Lambda_{\mathcal{A}}$, satisfies   $D_{H^{s}}(\Lambda_{\mathcal{A}}) = H^{1+s}(\R^d)$ holds with equivalent norms. 

\end{thm}

It is also shown that $\Lambda_{\mathcal{A}}$ admits a bounded $H^\infty$ calculus in $L^2 (\R^d)$; see Theorem \ref{thm.domain.DN}. 
Similar result is obtained in \cite{ES}, where they studied the 
Dirichlet-Neumann map for the Laplace operator in a bounded domain with
$C^{1+\alpha}$ boundary. 
See also \cite{Taylor2} for general properties
of the Dirichlet-Neumann map and the relation with 
the layer potentials.

This paper is organized as follows. In Section \ref{sec.preliminary} we state some general results on Poisson operators from \cite{MaekawaMiura1}, which plays a central role in our argument. Section \ref{sec.domain} is the core of this paper. In Section \ref{subsec.domain.L^2} we study the Poisson semigroup and its generator in $L^2 (\R^d)$ with the aid of the calculus of the symbol $\mu_{\mathcal{A}}$, while the Dirichlet-Neumann map in $L^2 (\R^d)$ is studied in Section \ref{subsec.domain.DN.L^2}. The analysis of these operators in $H^s(\R^d)$ is performed in Sections \ref{subsec.domain.H^1} - \ref{subsec.proof.thm}. As stated in Remark \ref{rem.pseudo}, our approach recovers some properties of the pseudo-differential operator $\mu_{\mathcal{A}}(\cdot,D_x)$ in $H^s(\R^d)$, which is stated in the appendix.

\section{Preliminaries}\label{sec.preliminary}

In this section we recall some results in \cite{MaekawaMiura1}. As stated in the introduction, the Poisson semigroup $\{E_{\mathcal{A}}(t)\}_{t\geq 0}$ defines a strongly continuous and analytic semigroup in $H^{1/2}(\R^d)$, and thus we have the representation $E_{\mathcal{A}}(t) = e^{-t\mathcal{P}_{\mathcal{A}}}$ with its generator $-\mathcal{P}_{\mathcal{A}}$. The next proposition gives the condition so that $\{e^{-t\mathcal{P}_{\mathcal{A}}}\}_{t\geq 0}$ is extended as a semigroup in $L^2 (\R^d)$.
\begin{prop}[{\cite[Proposition 3.3]{MaekawaMiura1}}]\label{prop.pre.1} The following two statements are equivalent.

\vspace{0.1cm}

\noindent {\rm (i)} $D_{H^{1/2}} (\mathcal{P}_{\mathcal{A}})\subset D_{L^2}(\Lambda_{\mathcal{A}^*})$ and $\| \Lambda_{\mathcal{A}^*} f\|_{L^2 (\R^d)} \leq C  \| f \|_{H^1 (\R^d)}$ holds for $f\in D_{H^{1/2}}(\mathcal{P}_{\mathcal{A}})$,

\noindent {\rm (ii)} $\{e^{-t\mathcal{P}_{\mathcal{A}}}\}_{t\geq 0}$ is extended as a strongly continuous semigroup in $L^2 (\R^d)$ and  $D_{L^2}(\mathcal{P}_{\mathcal{A}})$ is continuously embedded in $H^1 (\R^d)$. 

\vspace{0.1cm}

\noindent Moreover, if the condition {\rm (ii)} (and hence, {\rm (i)}) holds then $D_{L^2}(\mathcal{P}_{\mathcal{A}})$ is continuously embedded in $D_{L^2}(\Lambda_{\mathcal{A}})$,  $H^1 (\R^d)$ is continuously embedded in $D_{L^2}(\Lambda_{\mathcal{A}^*})$, and it follows that
\begin{align}
\mathcal{P}_{\mathcal{A}}  f   & = M_{1/b} \Lambda_{\mathcal{A}} f + M_{{\bf r_2}/b} \cdot \nabla_x f,\label{assume.prop.pre.1.1}\\
\langle A' \nabla_x f, \nabla_x g\rangle _{L^2(\R^d)} & = \langle \mathcal{P}_{\mathcal{A}} f, \Lambda_{\mathcal{A}^*} g + M_{\bf \bar{r}_1}\cdot \nabla_x g\rangle _{L^2(\R^d)}\label{assume.prop.pre.1.2}
\end{align}
for $f\in D_{L^2}(\mathcal{P}_{\mathcal{A}})$ and $g\in H^1(\R^d)$.

\end{prop}

In order to show $D_{L^2}(\mathcal{P}_{\mathcal{A}}) = H^1 (\R^d)$ we will use 
\begin{prop}[{\cite[Corollary 3.5, Proposition 3.6]{MaekawaMiura1}}]\label{prop.pre.2}  Assume that $\{e^{-t\mathcal{P}_{\mathcal{A}}}\}_{t\geq 0}$ and $\{e^{-t\mathcal{P}_{\mathcal{A}^*}}\}_{t\geq 0}$ are extended as  strongly continuous semigroups in $L^2 (\R^d)$ and that $D_{L^2}(\mathcal{P}_{\mathcal{A}})$ and $D_{L^2}(\mathcal{P}_{\mathcal{A}^*})$ are continuously embedded in $H^1 (\R^d)$. Then we have 
\begin{align}
& \langle A'\nabla_x f, \nabla_x g\rangle _{L^2 (\R^d)} = \langle \mathcal{P}_{\mathcal{A}}f, M_{\bar{b}} \mathcal{P}_{\mathcal{A}^*} g\rangle _{L^2 (\R^d)},~~~~~~~~ f\in D_{L^2}(\mathcal{P}_{\mathcal{A}}), ~~g\in D_{L^2}(\mathcal{P}_{\mathcal{A}^*}),\\
& C' \| f \|_{H^1 (\R^d)} \leq \| \mathcal{P}_{\mathcal{A}} f\|_{L^2 (\R^d)} + \| f\|_{L^2(\R^d)} \leq C  \| f\|_{H^1 (\R^d)},~~~~~~~~~f\in D_{L^2}(\mathcal{P}_{\mathcal{A}}).
\end{align}
If in addition that $\displaystyle \liminf_{t\rightarrow 0} \| \dd / \dd t ~ e^{-t\mathcal{P}_{\mathcal{A}}} f \|_{L^2 (\R^d)} <\infty$ holds for all $f\in C_0^\infty (\R^d)$ then $D_{L^2}(\mathcal{P}_{\mathcal{A}}) = H^1 (\R^d)$ with equivalent norms.

\end{prop}

As for the factorizations of $\mathcal{A}'$ and $\mathcal{A}$, we have 
\begin{prop}[{\cite[Lemma 3.7]{MaekawaMiura1}}]\label{prop.pre.3} Assume that the semigroups $\{ e^{-t\mathcal{P}_{\mathcal{A}}}\}_{t\geq 0}$  and $\{ e^{-t\mathcal{P}_{\mathcal{A}^*}}\}_{t\geq 0}$ in $H^{1/2}(\R^d)$ are extended as  strongly continuous semigroups in $L^2 (\R^d)$ and that $D_{L^2}(\mathcal{P}_{\mathcal{A}})= D_{L^2}(\mathcal{P}_{\mathcal{A}^*}) = H^1 (\R^d)$ holds with equivalent norms. Then $H^1 (\R^d)$ is continuously embedded in $D_{L^2}(\Lambda_{\mathcal{A}})\cap D_{L^2}(\Lambda_{\mathcal{A}^*})$ and 
\begin{align}
\mathcal{P}_{\mathcal{A}} f & = M_{1/b} \Lambda_{\mathcal{A}} f + M_{{\bf r_2}/b}\cdot \nabla_x f,~~~~~~~~~~f\in H^1 (\R^d),\label{eq.prop.pre.3.1} \\
\mathcal{P}_{\mathcal{A}^*} g &  = M_{1/\bar{b}} \Lambda_{\mathcal{A}^*} g + M_{{\bf \bar{r}_1}/\bar{b}}\cdot \nabla_x g,~~~~~~~~~~ g\in H^1 (\R^d).\label{eq.prop.pre.3.2}
\end{align}
Moreover, the realizations of $\mathcal{A}'$ in $L^2 (\R^d)$ and of $\mathcal{A}$ in $L^2 (\R^{d+1})$ are respectively factorized as
\begin{align}
\mathcal{A}' & = M_b \mathcal{Q}_{\mathcal{A}} \mathcal{P}_{\mathcal{A}},~~~~~~~\quad \mathcal{Q}_{\mathcal{A}} = M_{1/b} ( M_{\bar{b}} \mathcal{P}_{\mathcal{A}^*} )^*,\label{eq.prop.pre.3.3'} \\
\mathcal{A}  & = - M_b (\partial _t - \mathcal{Q}_{\mathcal{A}} ) ( \partial_t + \mathcal{P}_{\mathcal{A}}).\label{eq.prop.pre.3.3}
\end{align}
Here $( M_{\bar{b}} \mathcal{P}_{\mathcal{A}^*} )^*$ is the adjoint of $ M_{\bar{b}} \mathcal{P}_{\mathcal{A}^*}$ in $L^2 (\R^d)$.

\end{prop}

\section{Analysis of Poisson operator in $H^s (\R^d)$}\label{sec.domain}

To study the Poisson operator we consider the boundary value problem 
\begin{equation}\label{eq.dirichlet}
\begin{cases}
& \mathcal{A} u  = F~~~~~~~~{\rm in}~~~\R^{d+1}_+,\\
& ~~u  = g~~~~~~~~~{\rm on}~~\partial\R^{d+1}_+.
\end{cases}
\end{equation}
Let $x,\xi\in \R^d$ and let $\mu_{\mathcal{A}}=\mu_{\mathcal{A}}(x,\xi)\in \{\mu \in \C ~|~{\rm Im} \,\mu >0\}$ be the root of 
\begin{align}
b (x) \mu^2 + \big ({\bf r_1} (x) + {\bf r_2} (x) \big ) \cdot \xi \mu + \langle A'(x) \xi,\xi\rangle =0.\label{eq.root}
\end{align}   
Then we have
\begin{align}
\mu_{\mathcal{A}} (x,\xi) = - \frac{{\bf v} (x) \cdot \xi}{2} + i \big \{ \frac{1}{b (x) } \langle A'(x) \xi,\xi\rangle -\frac{1}{4}({\bf v}(x) \cdot \xi )^2\big \}^\frac12,~~~~~~~~{\bf v} = \frac{{\bf r_1} + {\bf r_2}}{b}.\label{def.mu_A}
\end{align}
Here the square root in \eqref{def.mu_A} is taken as the principal branch. From \eqref{ellipticity} one can check the estimates 
\begin{align}
|\mu_{\mathcal{A}} (x,\xi) |   \leq C |\xi |, ~~~&~~~ {\rm Im} \,\mu_{\mathcal{A}} (x,\xi) \geq C' |\xi |, \label{estimate.mu_A.1}\\
{\rm Re} \big ( \langle A'\xi,\xi\rangle - \frac{b}{4} ({\bf v}\cdot\xi)^2 \big ) & \geq \nu_1 \big ( |\xi|^2 + \frac{|{\bf v}\cdot \xi |^2}{4} \big ), \label{estimate.mu_A.2}
\end{align}
where $C, C'$ are positive constants depending only on $\nu_1, \nu_2$. As is well known, $\mu_{\mathcal{A}}$ describes the principal symbol of the Poisson operator.

\subsection{Domain of Poisson operator in $L^2 (\R^d)$}\label{subsec.domain.L^2}

The aim of this section is to prove that the domain of the Poisson operator in $L^2 (\R^d)$ is $H^1 (\R^d)$. For a given $h\in \mathcal{S}(\R^d)$ we set 
\begin{align}
\big ( U_{\mathcal{A},0} (t) h \big )(x) = \frac{1}{(2 \pi)^\frac{d}{2}} \int_{\R^d} e^{i t \mu_{\mathcal{A}} (x,\xi) + ix\cdot \xi } \hat{h}(\xi) \dd \xi,\label{def.U_A0}
\end{align}
where $\hat{h}$ is the Fourier transform of $h$. The operator $U_{\mathcal{A},0}(t)$ represents the principal part of the Poisson semigroup, and we first give some estimates of  $U_{\mathcal{A},0}(t)$. To this end let us introduce the operator 
\begin{align}
\big ( G_p (t) h \big ) (x) =  \frac{1}{(2\pi)^\frac{d}{2}}\int_{\R^d} p (x,\xi,t ) e^{it\mu_{\mathcal{A}} (x,\xi) + ix\cdot \xi} \hat{h} (\xi )\dd\xi,\label{def.G_p}
\end{align}
for a given measurable function $p=p(x,\xi,t)$ on $\R^d\times \R^d \times \R_+$. 
\begin{lem}\label{lem.G_p} Let $T\in (0,\infty]$. Assume that $p=p(x,\xi,t)$ satisfies 
\begin{align}
\sup_{x\in \R^d,\xi\ne 0}\sup_{0<t<T} ~ \big ( \sum_{k=0}^{d+1}  (1 + t  |\xi| )^{-l_k} |\xi |^{k} |\nabla_\xi^k p (x,\xi,t) | + (t  |\xi| )^{-l_{j_0}} |\xi |^{j_0} |\nabla_\xi^{j_0} p (x,\xi,t) | \big ) \leq L<\infty \label{assume.lem.G_p.1}
\end{align}
for some $l_k\geq 0$ and for some $l_{j_0}>0$, $j_0\in \{0,\cdots,d\}$. Then we have 
\begin{align}
\sup_{0<t<T} \| G_p (t)  h\|_{L^2 (\R^d)} \leq C \| h \|_{L^2 (\R^d)},\label{est.lem.G_p.1}
\end{align}
where $C$ depends only on $d$, $\nu_1$, $\nu_2$, $l_k$, $l_{j_0}$, and $L$.  Furthermore, if $p$ satisfies 
\begin{align}
\sup_{x\in \R^d, \xi\ne 0} \sup_{t>0} ~ \sum_{k=0}^{d+1}  (t  |\xi| )^{-l_k} |\xi |^{k} |\nabla_\xi^{k} p (x,\xi,t) | \leq L<\infty \label{assume.lem.G_p.2}
\end{align}
for some  $l_k>0$, $k=0,\cdots, d+1$, then  
\begin{align} 
\int_0^\infty \| G_p (t) h \|_{L^2 (\R^d)}^2 \frac{\dd t}{t} \leq C' \| h \|_{L^2 (\R^d)}^2,\label{est.lem.G_p.2}
\end{align}
where $C'$ depends only on $d$, $\nu_1$, $\nu_2$, $l_k$, and $L$.

\end{lem}
The proof of Lemma \ref{lem.G_p} is rather standard and will be stated in the appendix for convenience to the reader. Now we have 
\begin{lem}\label{lem.U_0}  Let $k\in \N\cup \{0\}$. For $h\in \mathcal{S} (\R^d)$ and $t>0$ it follows that 
\begin{align}
&  \| t^k \frac{\dd^k}{\dd t^k}U_{\mathcal{A}, 0} (t) h \|_{L^2 (\R^d)} +  \| t e^{-t} \nabla_x U_{\mathcal{A},0} (t) h \|_{L^2 (\R^d)} + \| t  U_{\mathcal{A},0} (t) (-\Delta_x)^\frac12 h \|_{L^2 (\R^d)} \nonumber \\
&{\phantom{ \| U_{\mathcal{A}, 0} (t) h \|_{L^2 (\R^d)} +  \| t \nabla U_{\mathcal{A},0} (t) h \|_{L^2 (\R^d)}}} ~~~~~~~~ + \| [U_{\mathcal{A},0} (t), \nabla_x] h \|_{L^2 (\R^d)}  \leq C \| h \|_{L^2 (\R^d)},\label{est.lem.U_0.1}
\end{align}
where $[B_1,B_2]$ is the commutator of the operators $B_1$, $B_2$, and 
\begin{align}
\int_0^\infty \| e^{-t} U_{\mathcal{A},0} (t) h \|_{\dot{H}^\frac12 (\R^d)}^2 \dd t \leq C \| h \|_{L^2 (\R^d)}^2. \label{est.lem.U_0.3}
\end{align} 
In particular, $\displaystyle \lim_{t\rightarrow 0} U_{\mathcal{A},0} (t) h =h$ in $L^2 (\R^d)$ for any $h\in L^2 (\R^d)$.
\end{lem}

\noindent {\it Proof.} The estimate \eqref{est.lem.U_0.1} is a direct consequence of Lemma \ref{lem.G_p}. For example, we take $p(x,\xi,t) =t  |\xi|$ for the estimate of $t U_{\mathcal{A},0} (t) (-\Delta_x)^{1/2} h$, and take $p (x,\xi,t) = i t \nabla_x \mu_{\mathcal{A}} (x,\xi)$ for $[U_{\mathcal{A},0}, \nabla_x ]h$, and so on.  As for \eqref{est.lem.U_0.3}, we use the Schur lemma as in the proof of \eqref{est.lem.G_p.2}.
By \eqref{est.lem.U_0.1} and $\| f \|_{\dot{H}^{1/2}}\leq \| f \|_{L^2}^{1/2} \| \nabla_x f \|_{L^2}^{1/2}$ it is easy to see that $t^{1/2} e^{-t} \| U_{\mathcal{A},0} (t) h \|_{\dot{H}^{1/2} (\R^d)} \leq C \| h \|_{L^2 (\R^d)}$ for all $t>0$. Let $t\geq s>0$, and let $\psi_s$ be the function defined in the proof of Lemma \ref{lem.G_p} in the appendix.  Then we have from $\psi_s = \Delta_x \tilde \psi_s$ and \eqref{est.lem.U_0.1},
\begin{align*}
t^\frac12 \| e^{-t} U_{\mathcal{A},0} (t) \psi_s * h \|_{\dot{H}^\frac12 (\R^d)} & \leq t^\frac12 \| e^{-t} U_{\mathcal{A},0} (t) \Delta_x \tilde \psi_s * h \|_{L^2 (\R^d)}^\frac12  \| \nabla_x e^{-t} U_{\mathcal{A},0} (t) \psi_s * h \|_{L^2 (\R^d)}^\frac12\\
&  \leq C t^{-\frac12} \| \nabla_x \tilde \psi_s * h \|_{L^2 (\R^d)}^\frac12 \| h \|_{L^2 (\R^d)}^\frac12 \leq C t^{-\frac12} s^\frac12 \| h \|_{L^2 (\R^d)}.
\end{align*}
Let $s\geq t>0$. Then the relation $\nabla_x U_{\mathcal{A},0} (t)  = t^{1/2} G_p (t) (-\Delta_x)^{1/4}  + U_{\mathcal{A},0} (t) \nabla_x$ with $p = i t^{1/2}|\xi|^{-1/2} \nabla_x \mu_{\mathcal{A}} $ combined with \eqref{est.lem.G_p.1} and  \eqref{est.lem.U_0.1} yields
\begin{align*}
\|  \nabla_x U_{\mathcal{A},0} (t) h \|_{L^2(\R^d)}\leq C t^\frac12 \| (-\Delta_x)^\frac14 \psi_s * h \|_{L^2(\R^d)} + C \| \nabla_x \psi_s * h \|_{L^2 (\R^d)} \leq C( t^\frac12 s^{-\frac12} + s^{-1} ),
\end{align*}
which implies
\begin{align*}
t^\frac12 \| e^{-t} U_{\mathcal{A},0} (t) \psi_s * h \|_{\dot{H}^\frac12 (\R^d)} \leq C t^\frac12 e^{-t} ( t^\frac14 s^{-\frac14} + s^{-\frac12} ) \| h \|_{L^2 (\R^d)} \leq C t^\frac14 s^{-\frac14} \| h \|_{L^2 (\R^d)}.
\end{align*}
Collecting these above, we can apply the Schur lemma \cite[pp.643-644]{Grafakos} to $\{ t^{1/2} e^{-t} (-\Delta_x)^{1/4} U_{\mathcal{A},0} (t)\}_{t>0}$ to obtain \eqref{est.lem.U_0.3}. The last statement of the proposition follows from \eqref{est.lem.U_0.1} and the density argument. The proof is complete. 

\vspace{0.5cm}

We look for a solution $u$ to \eqref{eq.dirichlet} with $F=0$ and $g=h$  of the form 
\begin{align}
u=M_\chi U_{\mathcal{A},0} h + U_{\mathcal{A},1}  h,
\end{align}
where $\chi=\chi(t)$ is a smooth cut-off function such that $\chi(t)=1$ if $t\in [0,1]$ and $\chi(t) =0$ if $t\geq 2$, and $U_{\mathcal{A},1} h$ is a solution to \eqref{eq.dirichlet} with $F= - \mathcal{A} (M_ \chi U_{\mathcal{A},0} h)$ and $g=0$.

\begin{lem}\label{lem.U_1} For any $h\in \mathcal{S} (\R^d)$ there exists a unique  solution $U_{\mathcal{A},1} h\in \dot{H}_0^1(\R^{d+1}_+)$  to \eqref{eq.dirichlet} with $F= - \mathcal{A} ( M_\chi U_{\mathcal{A},0} h)$ and $g=0$, which satisfies
\begin{align}
\| \nabla U_{\mathcal{A},1} h\|_{L^2 (\R^{d+1}_+)}\leq C \| (I - \Delta_x)^{-\frac14} h \|_{L^2 (\R^d)}.\label{est.lem.U_1.1}
\end{align}

\end{lem}

\noindent {\it Proof.} For simplicity we write $U_0$ and $U_1$ for $U_{\mathcal{A},0}$ and $U_{\mathcal{A},1}$. We set
$$A'=(a_{i,j})_{1\leq i,j\leq d},~~~~~~~~~~ ~{\bf a'} = \nabla_x\cdot A'= (\sum_{1\leq k\leq d}\partial_k a_{k,j})_{1\leq j\leq d},$$
and 
\begin{equation}
\Pi h  = 
\begin{pmatrix}
&  \Pi' h \mbox{} ~\\
& \Pi_{d+1} h  \mbox{} ~ 
\end{pmatrix}
=
\begin{pmatrix}
& \displaystyle G_{ it A' \nabla_x\mu_{\mathcal{A}} } h \\
& \displaystyle   G_{\nabla_x\cdot {\bf r_1}} h + G_{it ({\bf r_1} +{\bf r_2} )\cdot \nabla_x\mu_{\mathcal{A}}} h 
\end{pmatrix}
\in \C^{d+1}.\label{proof.lem.U_1.0}
\end{equation}
Then a direct computation yields
\begin{align*}
-\mathcal{A} U_0 h  = \nabla\cdot \Pi h +  G_\zeta h, ~~~~~~\zeta (x,\xi, t) & =   i \big ({\bf r_1} (x) + {\bf r_2}(x)  + i t A' (x) \xi \big ) \cdot \nabla_x \mu_{\mathcal{A}} (x,\xi) + i {\bf a'} (x) \cdot \xi.
\end{align*}
Hence $U_1h$ should be constructed as the solution of \eqref{eq.dirichlet} with $g=0$ and 
\begin{align}
F & = - \mathcal{A} (M_\chi U_0 h)  = \nabla \cdot M_\chi \Pi h + M_\chi G_\zeta h  + R h,\label{proof.lem.U_1.1} \\
R h & =  \nabla_x \cdot M_{({\bf r_1} + {\bf r_2})\partial_t \chi}  U_0  h+ \partial_t M_{b \partial_t\chi} U_0  h - M_{\partial_t\chi } \Pi_{d+1} h - M_{b\partial_t^2\chi + \nabla_x\cdot {\bf r_2}\partial_t\chi } U_0 h.\nonumber 
\end{align}
To obtain \eqref{est.lem.U_1.1}, 
let us estimate each term of $F$ in $\dot{H}^{-1}(\R^{d+1}_+)$. 
Note that $R h$ is supported in $\{(x,t)\in \R^{d+1}_+~|~1\leq t\leq 2\}$ by the definition of $\chi$. In particular, it is not difficult to show
\begin{align}
| \langle R h, \varphi \rangle _{L^2 (\R^{d+1}_+)} | \leq C \| (I-\Delta_x)^{-\frac14} h \|_{L^2 (\R^d)} \| \nabla \varphi \|_{L^2 (\R^{d+1}_+)},~~~~~~\varphi\in H_0^1(\R^{d+1}_+).\label{proof.lem.U_1.2}
\end{align}
Thus we focus on the leading  terms $\nabla \cdot M_\chi \Pi h$ and  $M_\chi G_\zeta h$. By using Lemma \ref{lem.G_p} one can easily check the estimates
\begin{align*}
|\langle \nabla \cdot M_\chi \Pi h , \varphi \rangle _{L^2 (\R^{d+1}_+)} |& \leq \| M_\chi \Pi h \|_{L^2 (\R^{d+1}_+)} \| \nabla \varphi \|_{L^2 (\R^{d+1}_+)} \leq C \| h \|_{L^2 (\R^d)} \| \nabla \varphi \|_{L^2 (\R^{d+1}_+)},\\ 
|\langle M_\chi G_\zeta h, \varphi \rangle _{L^2 (\R^{d+1}_+)} | &  =  | \int_0^\infty \langle M_\chi G_\zeta h, \int_0^t \partial_s \varphi \dd s \rangle _{L^2 (\R^d)} \dd t | \\
&  \leq \int_0^\infty t^\frac12 \| M_\chi G_\zeta h \|_{L^2 (\R^d)}\dd t \| \nabla \varphi \|_{L^2(\R^{d+1}_+)}\leq C \| h \|_{L^2 (\R^d)}  \| \nabla \varphi \|_{L^2 (\R^{d+1}_+)}, 
\end{align*}
for $\varphi \in H^1_0 (\R^{d+1}_+)$. 
Next we consider the estimate of $M_\chi G_\zeta (-\Delta_x)^\frac14 h$. 
By using the following general relation for $G_p$ that 
\begin{align}
G_p (t) h =\partial_t G_{p/(i\mu_{\mathcal{A}})} (t) h - G_{\partial_t p / (i\mu_{\mathcal{A}})} (t)h, ~~~~~~G_p (t) (-\Delta_x)^\frac14 h = G_{p |\xi|^{1/2}} (t) h,
\label{proof.lem.U_1.4}
\end{align}
and by using  $\partial_t^2 \zeta =0$ 
and $\partial_t\mu_{\mathcal{A}}=0$, 
we observe that it suffices to estimate 
\[
\partial_t M_\chi G_{\zeta |\xi|^{1/2} /(i\mu_{\mathcal{A}}) } h + \partial_t M_\chi G_{\partial_t \zeta |\xi|^{1/2}/\mu_{\mathcal{A}}^2} h,
\]
for the other terms are of lower order. We see 
\begin{align*}
|\langle \partial_t M_\chi G_{\zeta |\xi|^{1/2} /(i\mu_{\mathcal{A}}) } h + \partial_t M_\chi G_{\partial_t \zeta |\xi|^{1/2}/\mu_{\mathcal{A}}^2} h, \varphi \rangle _{L^2 (\R^{d+1}_+)} | & \leq C \| M_\chi   G_{\tilde \zeta } h\|_{L^2 (\R^{d+1}_+)} \| \nabla \varphi \|_{L^2 (\R^{d+1}_+)}
\end{align*}
for $\varphi\in H^1_0(\R^{d+1}_+)$, where $\tilde \zeta = \zeta |\xi|^{1/2} /(i\mu_{\mathcal{A}}) + \partial_t \zeta |\xi|^{1/2}/\mu_{\mathcal{A}}^2$. Then the bound of  $\| M_\chi   G_{\tilde \zeta } h\|_{L^2 (\R^{d+1}_+)}$ is reduced 
to that of 
$\int_{\R_+} \| G_p (t) h \|_{L^2(\R^d)}^2 t^{-1} \dd t$ with $p$ satisfying \eqref{assume.lem.G_p.2}, and hence, $\| M_\chi   G_{\tilde \zeta } h\|_{L^2 (\R^{d+1}_+)} \leq C \| h \|_{L^2 (\R^d)}$ follows 
from \eqref{est.lem.G_p.2}, as desired. Similarly, we have 
\begin{align*}
 \| M_\chi \Pi (-\Delta_x)^\frac14 h \|_{L^2 (\R^{d+1}_+)} \leq C \| h \|_{L^2 (\R^d)}.
\end{align*}
Collecting these above, we arrive at 
\begin{align}
| \langle - \mathcal{A} (M_ \chi U_0 h), \varphi\rangle _{L^2 (\R^{d+1}_+)} |\leq C \| (I-\Delta_x )^{-\frac14} h \|_{L^2 (\R^d)} \| \nabla\varphi \|_{L^2(\R^{d+1}_+)},~~~~~~~\varphi \in H^1_0 (\R^{d+1}_+).\label{proof.lem.U_1.5}
\end{align}
Thus, by the Lax-Milgram theorem there is a unique solution $U_1h\in \dot{H}^1_0 (\R^{d+1}_+)$ to \eqref{eq.dirichlet} with $F=- \mathcal{A} (M_\chi U_0 h)$ and $g=0$, which satisfies \eqref{est.lem.U_1.1}.   The proof is complete.

\vspace{0.5cm}

Lemmas \ref{lem.U_0} and \ref{lem.U_1} imply the estimate $\| e^{-t\mathcal{P}_{\mathcal{A}}} h \|_{L^2 (\R^d)} \leq C \| h\|_{L^2 (\R^d)}$ for $t\in (0,1]$ and $h\in H^1 (\R^d)$. Since the same argument can be applied for  $\{e^{-t\mathcal{P}_{\mathcal{A}^*}}\}_{t\geq 0}$, we have
\begin{cor}\label{cor.lem.U_0.U_1} The semigroups $\{e^{-t\mathcal{P}_{\mathcal{A}}}\}_{t\geq 0}$ and $\{e^{-t\mathcal{P}_{\mathcal{A}^*}}\}_{t\geq 0}$ in $H^{1/2}(\R^d)$ are extended as strongly continuous analytic semigroups in $L^2 (\R^d)$. Moreover, $H^1 (\R^d)$ is continuously embedded in $D_{L^2}(\Lambda_{\mathcal{A}})\cap D_{L^2}(\Lambda_{\mathcal{A}^*})$, and  $D_{L^2}(\mathcal{P}_{\mathcal{A}})$ and $D_{L^2}(\mathcal{P}_{\mathcal{A}^*})$ are continuously embedded in $H^1 (\R^d)$. We also have the estimate
\begin{align}
\| e^{-t\mathcal{P}_{\mathcal{A}}} f \|_{L^2 (\R^d)} \leq C t^{-\frac12} \| f \|_{H^{-\frac12} (\R^d)},~~~~~~~~~0<t<1.\label{est.cor.lem.U_0.U_1}
\end{align}

\end{cor}

\noindent {\it Proof.} By the variational characterization of $\dot{H}^{1/2}(\R^d)$, the estimate \eqref{est.lem.U_1.1} implies that 
\begin{align}
\| U_{\mathcal{A},1} (t) h \|_{\dot{H}^\frac12 (\R^d)} \leq \| \nabla U_{\mathcal{A},1} h \|_{L^2 (\R^{d+1}_+)} \leq C \| h \|_{H^{-\frac12} (\R^d)}, ~~~~~~~~~t>0.\label{proof.cor.lem.U_0.U_1.1}
\end{align}
Thus \eqref{est.lem.U_0.3} and \eqref{proof.cor.lem.U_0.U_1.1} verifies the condition (i) of \cite[Proposition 4.3]{MaekawaMiura1}, which gives $D_{L^2}(\mathcal{P}_{\mathcal{A}})\hookrightarrow H^1 (\R^d)$. The same is true for $D_{L^2}(\mathcal{P}_{\mathcal{A}^*})$, and then we also obtain the embedding $H^1 (\R^d)\hookrightarrow D_{L^2}(\Lambda_{\mathcal{A}}) \cap D_{L^2}(\Lambda_{\mathcal{A}^*})$ by Proposition \ref{prop.pre.1}. Now it remains to show \eqref{est.cor.lem.U_0.U_1}. Let us recall the representation $e^{-t\mathcal{P}_{\mathcal{A}}} f = M_\chi U_{\mathcal{A},0} (t) f + U_{\mathcal{A},1} (t)f$. By the definition of $G_p$ in \eqref{def.G_p}  we have $U_{\mathcal{A},0} (t) f = t^{-1/2} G_{t^{1/2}\langle \xi \rangle^{1/2}} (I - \Delta_x)^{-1/4}f$, where $\langle \xi \rangle = (1+|\xi|^2)^{1/2}$. Then it is easy to see that  $p(x,\xi,t)=t^{1/2}\langle \xi \rangle^{1/2}$ satisfies \eqref{assume.lem.G_p.1} with $T=1$,  and therefore,
\begin{align*}
\| U_{\mathcal{A},0} (t) f \|_{L^2 (\R^d)}  =  t^{-\frac12} \| G_{t^{1/2}\langle \xi \rangle^{1/2}} (I - \Delta_x)^{-1/4}f\|_{L^2 (\R^d)} \leq Ct^{-\frac12} \| (I - \Delta_x)^{-1/4}f\|_{L^2 (\R^d)}.
\end{align*}
On the other hand, we have already proved the desired estimate for $U_{\mathcal{A},1} (t) f$ by \eqref{est.lem.U_1.1}. The proof is complete.

\vspace{0.5cm}

In order to establish the characterization of $D_{L^2}(\mathcal{P}_{\mathcal{A}})$ and  $D_{L^2}(\mathcal{P}_{\mathcal{A}^*})$, 
we need further estimates of $U_{\mathcal{A},1}$ as follows.

\begin{lem}\label{lem.U_1'} For any $h\in H^{1/2}(\R^d)$ we have $\dd /\dd t ~ U_{\mathcal{A},1} (t) h\in C ([0,\infty); H^{1/2} (\R^d))$ and 
\begin{align}
\sup_{t>0} \| \frac{\dd }{\dd t} U_{\mathcal{A},1} (t) h \|_{H^\frac12 (\R^d)} \leq C \| h \|_{H^\frac12 (\R^d)}. \label{est.lem.U_1'.1}
\end{align}

\end{lem}

\noindent {\it Proof.} As in the proof of Lemma \ref{lem.U_1} we write $U_0$ and $U_1$ for $U_{\mathcal{A},0}$ and $U_{\mathcal{A},1}$. First we assume that  $h\in D_{L^2}(\mathcal{P}_{\mathcal{A}})$. Let us recall that  $U_1 (t)h\in \dot{H}^1_0(\R^{d+1}_+)$ solves \eqref{eq.dirichlet} with $F= - \mathcal{A}(M_\chi U_0 h) = \nabla\cdot M_\chi \Pi h + M_\chi G_\zeta h + R h$ and $g=0$. 
We first see $\lim_{\delta\rightarrow 0}d/dt U_1(\delta)h$ exists in 
$H^{-1}(\R^d)$. Set $\R^{d+1}_{\delta,+}=\{(x,t)\in \R^{d+1}~|~t>\delta\}$ for $\delta>0$. Since 
$\dd /\dd t~ U_1 h \in \cap_{\delta>0} H^1 (\R^{d+1}_{\delta,+})$, 
we have 
\begin{align}
\langle  \frac{\dd}{\dd t} U_1 (\delta ) h + M_{{\bf r_2}/b}\cdot \nabla_x  U_1 (\delta ) h, \phi \rangle _{L^2 (\R^{d})} & = - \langle A \nabla U_1 h , \nabla\varphi \rangle _{L^2 (\R^{d+1}_{\delta,+})}- \langle  M_\chi \Pi' h , \nabla_x \varphi \rangle _{L^2 (\R^{d+1}_{\delta,+})}\nonumber \\
& ~~~ + \langle \partial_t M_\chi \Pi_{d+1} h + M_\chi G_\zeta h + R h,\varphi \rangle  _{L^2 (\R^{d+1}_{\delta,+})}\label{proof.lem.U_1'.1}
\end{align}
for any $\delta\in (0,1/2)$ and $\varphi\in H^1 (\R^{d+1}_{\delta,+})$ with $\varphi|_{t=\delta} = M_{1/\bar{b}}\phi$, $\phi\in H^{1/2}(\R^d)$. By taking $\varphi = e^{-(t-\delta)\mathcal{P}_{\mathcal{A}^*}} M_{1/\bar{b}}\phi$, the similar calculation as in the proof of Lemma \ref{lem.U_1} yields the estimate for the right-hand side of \eqref{proof.lem.U_1'.1} from above by $C \| h \|_{H^{1/2}(\R^d)}\|M_{1/\bar{b}} \phi \|_{H^{1/2}(\R^d)}$, and hence by $C \| h \|_{H^{1/2}(\R^d)}\| \phi \|_{H^{1/2}(\R^d)}$. On the other hand, we have 
\begin{align}
|\langle M_{{\bf r_2}/b}\cdot \nabla_x  U_1 (\delta ) h, \phi \rangle _{L^2 (\R^{d})}| & \leq C \| \nabla_x U_1 (\delta ) h \|_{L^2 (\R^d)} \| \phi \|_{L^2 (\R^d)} \nonumber  \\
& = C \| \nabla_x \big ( e^{-\delta \mathcal{P}_{\mathcal{A}}} - U_0 (\delta ) \big ) h\|_{L^2 (\R^d)}  \| \phi \|_{L^2 (\R^d)} \nonumber \\
& \leq C \big ( \| [\nabla_x, e^{-\delta \mathcal{P}_{\mathcal{A}}} ] h \|_{L^2 (\R^d)} + \| [\nabla_x, U_0 (\delta ) ]  h\|_{L^2 (\R^d)} \big )   \| \phi \|_{L^2 (\R^d)}\nonumber \\
& ~~~ + C \|  U_1 (\delta )  \nabla_x h\|_{L^2 (\R^d)}  \| \phi \|_{L^2 (\R^d)}. \label{proof.lem.U_1'.2}
\end{align}
In particular, we see $\|  M_{{\bf r_2}/b}\cdot \nabla_x  U_1 (\delta ) h \|_{L^2(\R^d)}\rightarrow 0$ as $\delta\rightarrow 0$, for the facts $h\in D_{L^2}(\mathcal{P}_{\mathcal{A}})$ and $D_{L^2}(\mathcal{P}_{\mathcal{A}})\hookrightarrow H^1 (\R^d)$ by Corollary \ref{cor.lem.U_0.U_1} imply  $\displaystyle \lim_{\delta\rightarrow 0} \| [\nabla_x, e^{-\delta \mathcal{P}_{\mathcal{A}}} ] h \|_{L^2 (\R^d)}=0$, while the convergence $\displaystyle \lim_{\delta\rightarrow 0} \| [\nabla_x, U_0 (\delta ) ]  h\|_{L^2 (\R^d)}=0$ follows from \eqref{est.lem.U_0.1} and the density argument. The limit $\displaystyle \lim_{\delta\rightarrow 0} \|  U_1 (\delta )  \nabla_x h\|_{L^2 (\R^d)}=0$ follows from Lemma \ref{lem.U_1} with $h$ replaced by $\nabla_x h$. By applying \eqref{proof.lem.U_1'.1} to another $\delta'\in (0,1/2)$ and then by estimating the difference $\langle  \dd/\dd t ~ U_1 (\delta ) h -   \dd/\dd t ~  U_1 (\delta' ) h, \phi \rangle _{L^2 (\R^{d})}$ we conclude that $\dd/\dd t~ U_1 (\delta ) h$ converges to some limit, denoted by $S_{\mathcal{A},1} h$, in $H^{-1/2}(\R^d)$ as $\delta\rightarrow 0$. We claim that
\begin{align}
\| S_{\mathcal{A},1} h \|_{H^s(\R^d)}  \leq C_s \| h \|_{H^s (\R^d)},~~~~~~~~~h\in D_{L^2}(\mathcal{P}_{\mathcal{A}}),~~~~s\in (0,\frac12]. 
\label{proof.lem.U_1'.5}
\end{align}
To this end, first note that 
the following equality with the choice of 
$\varphi (t) = e^{-t\mathcal{P}_{\mathcal{A}^*}}  \phi$, 
$\phi\in H^{1/2}(\R^d)$ holds:
\begin{align}
\langle S_{\mathcal{A},1} h, M_{\bar{b}} \phi \rangle _{H^{-\frac12},H^\frac12} & = \lim_{\delta\rightarrow 0} \langle  \frac{\dd}{\dd t}  U_1 (\delta ) h, M_{\bar{b}} \phi \rangle _{H^{-\frac12},H^\frac12}  =\lim_{\delta\rightarrow 0} \langle   \frac{\dd}{\dd t} U_1 (\delta ) h, M_{\bar{b}} \varphi (\delta) \rangle _{H^{-\frac12},H^\frac12} \nonumber \\
& = \lim_{\delta\rightarrow 0} \langle M_b  \frac{\dd}{\dd t}  U_1 (\delta ) h + M_{{\bf r_2}}\cdot \nabla_x  U_1 (\delta ) h , \varphi (\delta) \rangle _{L^2 (\R^d )} \nonumber \\
&  = - \langle A \nabla U_1 h , \nabla\varphi \rangle _{L^2 (\R^{d+1}_{+})}- \langle  M_\chi \Pi' h , \nabla_x \varphi \rangle _{L^2 (\R^{d+1}_{+})}\nonumber \\
& ~~~ + \langle \partial_t M_\chi \Pi_{d+1} h + M_\chi G_\zeta h + R h,\varphi \rangle  _{L^2 (\R^{d+1}_{+})}\nonumber \\
& =- \langle  M_\chi \Pi' h , \nabla_x \varphi \rangle _{L^2 (\R^{d+1}_{+})} + \langle \partial_t M_\chi \Pi_{d+1} h + M_\chi G_\zeta h + R h,\varphi \rangle  _{L^2 (\R^{d+1}_{+})}. \label{proof.lem.U_1'.3}
\end{align}
Here for each line we have respectively used $\varphi (\delta) = e^{-\delta\mathcal{P}_{\mathcal{A}^*}}  \phi \rightarrow  \phi$ in $H^{1/2} (\R^d)$ for $\phi \in H^{1/2} (\R^d)$, $M_{{\bf r_2}/b}\cdot \nabla_x  U_1 (\delta ) h\rightarrow 0$ in $L^2 (\R^d)$, \eqref{proof.lem.U_1'.1}, and then the fact that $ \langle A \nabla U_1 h , \nabla\varphi \rangle _{L^2 (\R^{d+1}_{+})} =0$, since $U_1 h\in \dot{H}^1_0(\R^{d+1}_+)$ and $\varphi (t) = e^{-t\mathcal{P}_{\mathcal{A}}^*} \phi$. Hence our next task 
is to estimate the right-hand side of \eqref{proof.lem.U_1'.3}. Let $s\in (0,1/2]$. We will show 
\begin{align}
& ~~~ |\langle  M_\chi \Pi' h , \nabla_x \varphi \rangle _{L^2 (\R^{d+1}_{+})}| +  |\langle  \partial_t M_\chi G_{it ({\bf r_1} +{\bf r_2} )\cdot \nabla_x\mu_{\mathcal{A}}} h,  \varphi \rangle _{L^2 (\R^{d+1}_{+})}| \nonumber \\
& ~~~~~~~+  |\langle  M_\chi\partial_t   G_{\nabla_x\cdot {\bf r_2}} h +   M_\chi G_\zeta h,\varphi \rangle _{L^2 (\R^{d+1}_+)}| \leq C \| h \|_{H^s (\R^d)} \| (I-\Delta_x)^{-\frac{s}{2}} \phi \|_{L^2 (\R^d)}, \label{proof.lem.U_1'.4}
\end{align}
which gives \eqref{proof.lem.U_1'.5}, since the other terms in \eqref{proof.lem.U_1'.3} are of lower order and  easy to handle. We observe that $\varphi (t) = e^{-t\mathcal{P}_{\mathcal{A}^*}} \phi = M_\chi U_{\mathcal{A}^*,0} (t)  \phi + U_{\mathcal{A}^*,1}(t) \phi $, and thus, it suffices to check \eqref{proof.lem.U_1'.3}  with $\varphi$ replaced by $U_{\mathcal{A}^*,0}(t) \phi$, for $U_{\mathcal{A}^*,1}(t) \phi$ is of lower order if $s$ is less than or equal to $1/2$ thanks to Lemma \ref{lem.U_1}.  The argument below is based on the quadratic estimate  as in \eqref{est.lem.G_p.2}. Note that the counterpart of Lemma \ref{lem.G_p} is valid for $U_{\mathcal{A}^*,0}(t)$.  Firstly we see from the definition of $\Pi'$ in \eqref{proof.lem.U_1.0},
\begin{align*}
& ~~~ |\langle M_\chi \Pi' h , \nabla_x U_{\mathcal{A}^*,0} (t) (-\Delta_x)^\frac{s}{2} \phi \rangle _{L^2 (\R^{d+1}_{+})}| \\
& \leq \int_0^\infty M_\chi \| G_{it A' \nabla_x\mu_{\mathcal{A}}} (t) h \|_{L^2 (\R^d)} \| \nabla_x   U_{\mathcal{A}^*,0} (t) (-\Delta_x)^\frac{s}{2}\phi \|_{L^2 (\R^d)} \dd t \\
& \leq \big ( \int_0^\infty  M_\chi \| t^{1-s} G_{i A' \nabla_x\mu_{\mathcal{A}}} (t)  h \|_{L^2 (\R^d)}^2 \frac{\dd t}{t} \big )^\frac12 \big ( \int_0^\infty M_\chi \| t^{1+ s} \nabla_x  U_{\mathcal{A}^*,0} (t) (-\Delta_x)^\frac{s}{2} \phi \|_{L^2 (\R^d)}^2 \frac{\dd t}{t} \big )^\frac12\\
& \leq C \| h \|_{\dot{H}^s (\R^d)} \| \phi \|_{L^2 (\R^d)}, ~~~~~~s\in [0,\frac12 ],
\end{align*}
by applying Lemma \ref{lem.G_p}. Similarly, we have 
\begin{align*}
 & ~~~ |\langle  \partial_t M_\chi G_{it ({\bf r_1} +{\bf r_2} )\cdot \nabla_x\mu_{\mathcal{A}}} h,  \varphi \rangle _{L^2 (\R^{d+1}_{+})}|\nonumber\\ 
& = |\langle M_\chi G_{it ({\bf r_1} +{\bf r_2} )\cdot \nabla_x\mu_{\mathcal{A}}} h, \frac{\dd}{\dd t} U_{\mathcal{A}^*,0} (t) (-\Delta_x)^\frac{s}{2}\phi \rangle _{L^2 (\R^{d+1}_{+})}|  \leq C \| h \|_{\dot{H}^s (\R^d)} \| \phi \|_{L^2 (\R^d)}, ~~~~~~~s\in [0,\frac12 ].
\end{align*}
Finally we consider the term $M_\chi\partial_t   G_{\nabla_x\cdot {\bf r_2}} h +   M_\chi G_\zeta h = M_\chi G_{\eta} h$, $\eta = \mu_{\mathcal{A}} \nabla_x\cdot {\bf r_2} + \zeta$.  Lemma \ref{lem.G_p}  yields
\begin{align*}
B_s (h) : = \big ( \int_0^\infty M_\chi \| t^{1-s}  G_{\eta} (t)  h \|_{L^2(\R^d)}^2 \frac{\dd t}{t}\big )^\frac12  \leq C \| h \|_{\dot{H}^s (\R^d )},~~~~~~s\in [0,\frac12 ].
\end{align*} 
Thus it follows that 
\begin{align}
& |\langle  M_\chi G_{\eta} h ,   U_{\mathcal{A}^*,0} (t) \phi \rangle _{L^2 (\R^{d+1}_+)}|  \leq \int_0^\infty  M_\chi  \| G_{\eta}  (t) h \|_{L^2 (\R^d )}   \|  U_{\mathcal{A}^*,0} (t) \phi \|_{L^2 (\R^d)} \dd t \nonumber \\
&~~~~~~~~~~~~~~ \leq  \sup_{t>0}  \|  U_{\mathcal{A}^*,0} (t) (I-\Delta_x)^{-1} \phi \|_{L^2 (\R^d )} \int_0^\infty  M_\chi  \| G_{\eta}  (t) h \|_{L^2 (\R^d )} \dd t\nonumber  \\
& ~~~~~~~~~~~~~~~~~~ + \int_0^\infty  M_\chi  \| G_{\eta}  (t) h \|_{L^2 (\R^d )}   \|  U_{\mathcal{A}^*,0} (t) (-\Delta_x ) (I-\Delta_x)^{-1} \phi \|_{L^2 (\R^d)} \dd t,\nonumber
\end{align}
and the first term of the right-hand side is bounded from above by 
\[
C \|(I-\Delta_x)^{-1} \phi \|_{L^2 (\R^d )} \big (\int_0^\infty \chi (t) t^{2 s}  \frac{\dd t}{t} \big )^{\frac12} B_s (h)\leq C_s B_s (h) \|(I-\Delta_x)^{-1} \phi \|_{L^2 (\R^d )},
\]
while the second term is estimated by 
\begin{align*}
& ~~~B_s (h)   \big ( \int_0^\infty M_\chi \| t^{s} U_{\mathcal{A}^*,0} (t) (-\Delta_x)^{\frac{s}{2}}(-\Delta_x)^{1-\frac{s}{2}} (I-\Delta_x)^{-1} \phi \|_{L^2 (\R^d)}^2 \frac{\dd t}{t}\big )^\frac12 \\
& \leq C B_s (h) \| (-\Delta_x)^{1-\frac{s}{2}} (I-\Delta_x)^{-1} \phi \|_{L^2 (\R^d)}.
\end{align*}
Thus we have 
\begin{align}
|\langle  M_\chi G_{\eta} h ,   U_{\mathcal{A}^*,0} (t) \phi \rangle _{L^2 (\R^{d+1}_+)}|& \leq C_s    \|(I-\Delta_x)^{-\frac{s}{2}} \phi \|_{L^2 (\R^d )}  \| h \|_{\dot{H}^s (\R^d)},\nonumber 
\end{align} 
where $C_s$ is a constant which tends to $\infty$ as $s\rightarrow 0$. Collecting these above, we arrive at \eqref{proof.lem.U_1'.5}.

Now let $V_1 (t) h\in C([0,\infty); H^{1/2} (\R^d))$ be the  unique weak solution to \eqref{eq.dirichlet} with $F=- \partial_t\mathcal{A}(M_\chi U_0 h) = \nabla\cdot \partial_t  M_\chi  \Pi h + \partial_t M_\chi  G_\zeta h + \partial_t R h $ and $g=S_{\mathcal{A},1} h$. Hence it has the form 
\begin{align}
V_1 (t) h = e^{-t\mathcal{P}_{\mathcal{A}}} S_{\mathcal{A},1} h + W_{\mathcal{A}} (t) h = \chi U_0 (t) S_{\mathcal{A},1} h + U_1  (t) S_{\mathcal{A},1} h + W_{\mathcal{A}} (t) h,\label{proof.lem.U_1'.6}
\end{align}
where $W_{\mathcal{A}}  h\in \dot{H}^1_0 (\R^{d+1}_+)$ is a weak solution to \eqref{eq.dirichlet} with $F=  \nabla\cdot \partial_t  M_\chi  \Pi h + \partial_t M_\chi  G_\zeta h + \partial_t R h $ and $g=0$. Note that 
\begin{align}
\int_0^\infty M_\chi\big (  \| t^\frac12 \partial_t  \frac{\dd}{\dd t} \Pi   h   \|_{L^2 (\R^d)}^2 + \| t^\frac12 G_\zeta (t)  h \|_{L^2 (\R^d)}^2 \big )\frac{\dd t}{t} \leq C \| h \|_{\dot{H}^\frac12 (\R^d)}^2\nonumber
\end{align}
by Lemma \ref{lem.G_p}, which implies the estimate $\| \nabla W_{\mathcal{A}}  h \|_{L^2 (\R^{d+1}_+)}\leq C \| h \|_{\dot{H}^{1/2} (\R^d)}$, and thus,
\begin{align}
\| W_{\mathcal{A}} (t) h  \|_{\dot{H}^\frac12(\R^d)}\leq \| \nabla w \|_{L^2 (\R^{d+1}_{t,+})}\leq C \| h \|_{\dot{H}^{1/2} (\R^d)},~~~~~~~~t>0.\label{proof.lem.U_1'.7}
\end{align}
On the other hand, since $\{e^{-t\mathcal{P}_{\mathcal{A}}}\}_{t\geq 0}$ defines a strongly continuous semigroup in $H^{1/2}(\R^d)$, we have from  \eqref{proof.lem.U_1'.5},
\begin{align}
\| e^{-t\mathcal{P}_{\mathcal{A}}} S_{\mathcal{A},1} h \|_{H^\frac12 (\R^d)} \leq C \| S_{\mathcal{A},1}  h  \|_{H^\frac12 (\R^d)} \leq C \| h \|_{H^\frac12 (\R^d)},~~~~~~0<t\leq 2. \label{proof.lem.U_1'.8}
\end{align}
Finally let us prove $V_1 (t) h = \dd / \dd t~U_1 (t) h$. To see this, we note that, for each $\delta\in (0,1/2)$, $\dd / \dd t~U_1 (t+\delta) h $ is the unique weak solution to \eqref{eq.dirichlet} with $F = \tau_{\delta}  \big ( \nabla\cdot \partial_t  M_\chi  \Pi h + \partial_t M_\chi  G_\zeta h + \partial_t R h \big )$ and $g= \dd / \dd t~U_1 (\delta) h$, where $\tau_{\delta} f (t) =f(t+\delta)$. Hence we can write 
\begin{align}
\frac{\dd }{\dd t} U_1 (t+\delta )h = e^{-t\mathcal{P}_{\mathcal{A}}}   \frac{\dd }{ \dd t} U_1 (\delta) h + W_{\mathcal{A}}^{(\delta)} (t) h,~~~~~~~~\delta\in (0,1/2), \label{proof.lem.U_1'.9}
\end{align}
where $W_{\mathcal{A}}^{(\delta)} h$ is the solution to \eqref{eq.dirichlet} with $F = \tau_{\delta}  \big ( \nabla\cdot \partial_t  M_\chi \Pi h + \partial_t M_\chi  G_\zeta h + \partial_t R h \big )$ and $g=0$. Then, arguing as in the estimate of $W_{\mathcal{A}}h$, we obtain a uniform bound of $W_{\mathcal{A}}^{(\delta)}h$ in $\dot{H}_0^1(\R^{d+1}_+)$, and it is not difficult to show  $W_{\mathcal{A}}^{(\delta)}h$ weakly converges to $W_{\mathcal{A}} h$ in $\dot{H}_0^1 (\R^{d+1}_+)$. On the other hand, we have from \eqref{est.cor.lem.U_0.U_1},
\begin{align*}
\| e^{-t\mathcal{P}_{\mathcal{A}}} S_{\mathcal{A},1} h - e^{-t\mathcal{P}_{\mathcal{A}}}   \frac{\dd }{ \dd t} U_1 (\delta) h \|_{L^2 (\R^d)} \leq C t^{-\frac12} \| S_{\mathcal{A},1} h -   \frac{\dd }{ \dd t} U_1 (\delta) h \|_{H^{-\frac12}(\R^d)}, 
\end{align*}
which converges to zero as $\delta \rightarrow 0$ by the definition of $S_{\mathcal{A},1}h$. Hence, the right-hand side of \eqref{proof.lem.U_1'.9} converges to $V_1(t)h$ in $L^2_{loc} (\R^{d+1}_+)$.  On the other hand, since  $\dd /\dd t~ U_1 h \in \cap_{\delta>0} H^1 (\R^{d+1}_{\delta,+})$, the left-hand side of \eqref{proof.lem.U_1'.9} converges to $\dd /\dd t~ U_1 (t)h$ in $L^2 (\R^d)$ as $\delta\rightarrow 0$ for each $t>0$.  Thus we have $\dd / \dd t~U_1 (t) h= V_1 (t)h$, and then \eqref{proof.lem.U_1'.7} - \eqref{proof.lem.U_1'.8} implies 
$\| \dd/\dd t~ U_1 (t) h \|_{H^{1/2}(\R^d)}\leq C \| h \|_{H^{1/2} (\R^d)}$ for $0<t\leq 2$. For $t>2$ we have from the equality $\dd /\dd t ~ U_1 (t) h = \dd / \dd t~ e^{-t\mathcal{P}_{\mathcal{A}}} h$ that 
\begin{align*}
\| \frac{\dd}{\dd t} ~ U_1 (t) h\|_{H^\frac12 (\R^d)} \leq C \|  e^{-s\mathcal{P}_{\mathcal{A}}} h \mid _{s=1} \|_{L^2 (\R^d)} \leq C \| h \|_{L^2 (\R^d)}. 
\end{align*}
Hence, \eqref{est.lem.U_1'.1} holds when $h\in D_{L^2}(\mathcal{P}_{\mathcal{A}})$.  Note that $D_{L^2}(\mathcal{P}_{\mathcal{A}})$ is dense in $H^{1/2}(\R^d)$. Hence the estimate \eqref{est.lem.U_1'.1} for $h\in D_{L^2}(\mathcal{P}_{\mathcal{A}})$ is extended to all $h\in H^{1/2}(\R^d)$ by the density argument.  The proof is complete.

\begin{cor}\label{cor.lem.U_1'}  It follows that  $D_{L^2}(\mathcal{P}_{\mathcal{A}}) = D_{L^2}(\mathcal{P}_{\mathcal{A}^*}) = H^1(\R^d)$ with equivalent norms.

\end{cor}

\noindent {\it Proof.} Note that the same result as in Lemma \ref{lem.U_1'} is valid  for  $\dd /\dd t ~ U_{\mathcal{A}^*,1} (t)$. On the other hand, by Lemma \ref{lem.U_0} and the definition \eqref{def.U_A0}, it is straightforward to see 
\begin{align*}
\sup_{t>0} \|  \frac{\dd}{\dd t} U_{\mathcal{A},0}(t) h \|_{L^2 (\R^d)}  +\sup_{t>0}  \| \frac{\dd}{\dd t} U_{\mathcal{A}^*,0} (t) h \|_{L^2 (\R^d)}<\infty ~~~~~~~{\rm for}~~h\in C_0^\infty (\R^d)
\end{align*}
by the definition \eqref{def.U_A0}. Thus we have
\begin{align}
\sup_{t>0} \| \frac{\dd}{\dd t} e^{-t\mathcal{P}_{\mathcal{A}}} h \|_{L^2 (\R^d)}  +\sup_{t>0}  \| \frac{\dd}{\dd t} e^{-t\mathcal{P}_{\mathcal{A}^*}} h \|_{L^2 (\R^d)}<\infty ~~~~~~{\rm for}~~~ h\in C_0^\infty  (\R^d ). \label{proof.cor.lem.U_1'.1} 
\end{align}
Hence Proposition \ref{prop.pre.2} (ii), together with Corollary \ref{cor.lem.U_0.U_1}, shows  $D_{L^2}(\mathcal{P}_{\mathcal{A}}) = D_{L^2}(\mathcal{P}_{\mathcal{A}^*}) = H^1(\R^d)$ with equivalent norms. The proof is complete.

\vspace{0.5cm}

\begin{rem}\label{rem.pseudo}{\rm Let $\Phi:A \longmapsto \mu_{\mathcal{A}}$ be the map defined by \eqref{def.appendix.Phi}. The arguments of the present section essentially rely on the integration by parts technique, and in particular, we did not use the mapping properties of the pseudo-differential operator $\Phi (A)(\cdot, D_x)$ such as the equivalence $\| \Phi (A)(\cdot,D_x) f \|_{L^2(\R^d)} + \| f\|_{L^2(\R^d)} \simeq \| f \|_{H^1 (\R^d)}$.  Since the above proof implies the identity $-\mathcal{P}_{\mathcal{A}}=i\Phi (A) (\cdot,D_x) + S_{\mathcal{A},1}$, where $S_{\mathcal{A},1}=\displaystyle \lim_{t\rightarrow 0}\dd/\dd t~U_{\mathcal{A},1}(t)$ is a lower order operator, our result actually gives an alternative proof (although it is lengthy) of the mapping properties of  $\Phi (A) (\cdot,D_x)$ in $H^s(\R^d)$ which are well known in the theory of pseudo-differential operators with nonsmooth coefficients; cf. \cite{KumanogoNagase, Marschall, Taylor, Abels, ES}. Especially, the fact that $i\Phi (A) (\cdot,D_x)$ in $L^2 (\R^d)$ with the domain $H^1 (\R^d)$ generates a strongly continuous and analytic semigroup in $L^2 (\R^d)$ is recovered by regarding  $i\Phi (A) (\cdot,D_x)$  as a perturbation from $-\mathcal{P}_{\mathcal{A}}$. More precise statements will be given in Theorem \ref{thm.class.Phi}.

}
\end{rem}

\subsection{Domain of Dirichlet-Neumann map in $L^2 (\R^d)$}\label{subsec.domain.DN.L^2}

In this section we consider the domain of the Dirichlet-Neumann map in $L^2 (\R^d)$. The result is stated as follows.
\begin{thm}\label{thm.domain.DN} It follows that $D_{L^2}(\Lambda_{\mathcal{A}}) = D_{L^2} (\Lambda_{\mathcal{A}^*}) = H^1 (\R^d)$ with equivalent norms. Moreover, $\Lambda_{\mathcal{A}}$ (and hence, $\Lambda_{\mathcal{A}^*}$) admits a bounded $H^\infty$ calculus in $L^2 (\R^d)$.

\end{thm}

\noindent {\it Proof.}  It suffices to show $D_{L^2}(\Lambda_{\mathcal{A}}) = H^1 (\R^d)$. Let $\mu_{\mathcal{A}}$ be as in \eqref{def.mu_A}. Let $f\in H^1 (\R^d)$. Since we have already shown $D_{L^2}(\mathcal{P}_{\mathcal{A}}) =H^1 (\R^d)$ in the previous section, Proposition \ref{prop.pre.3} gives
\begin{align}
\mathcal{P}_{\mathcal{A}} f & = M_{1/b} \Lambda_{\mathcal{A}} f + M_{{\bf r_2}/b} \cdot \nabla_x f,\label{proof.thm.domain.DN.1}
\end{align}
while, as stated in Remark \ref{rem.pseudo}, we have
\begin{align}
\mathcal{P}_{\mathcal{A}} f & = - i \mu_{\mathcal{A}} (\cdot,D_x) f - S_{\mathcal{A},1} f = - i M_{1/b} \lambda_{\mathcal{A}} (\cdot,D_x) f - S_{\mathcal{A},1} f  + M_{{\bf r_2}/b} \cdot \nabla_x f,\label{proof.thm.domain.DN.2}
\end{align}
where $S_{\mathcal{A},1}$ is the linear operator given in the proof of Lemma \ref{lem.U_1'} and $\lambda_{\mathcal{A}} (x,\xi )= b(x) \mu_{\mathcal{A}} (x,\xi) + {\bf r_2} (x) \cdot \xi$. Now let us define the linear operator  $\mathcal{J}_{\mathcal{A}}$ in $L^2 (\R^d)$  by 
\begin{align}
D_{L^2} (\mathcal{J}_{\mathcal{A}}) = H^1 (\R^d), ~~~~~~~~~~ \mathcal{J}_{\mathcal{A}} f =  - i \lambda_{\mathcal{A}} (\cdot,D_x) f - M_b  S_{\mathcal{A},1} f,\label{proof.thm.domain.DN.3}
\end{align} 
which gives $\Lambda_{\mathcal{A}} f = \mathcal{J}_{\mathcal{A}}f $ for $f\in H^1 (\R^d)$ by \eqref{proof.thm.domain.DN.1} - \eqref{proof.thm.domain.DN.2}. Thanks to Lemma \ref{lem.invariant.lambda} together with Remark \ref{rem.pseudo} the operator $i \lambda_{\mathcal{A}} (\cdot,D_x)$ in $L^2 (\R^d)$ with $D_{L^2}(\lambda_{\mathcal{A}} (\cdot,D_x)) =H^1 (\R^d)$ generates a strongly continuous and analytic semigroup in $L^2 (\R^d)$, and $-i \lambda_{\mathcal{A}} (\cdot,D_x)$ admits a bounded $H^\infty$ calculus in $L^2 (\R^d)$. On the other hand, we have from \eqref{proof.lem.U_1'.5}, 
\begin{align*}
\| M_b  S_{\mathcal{A},1} f\|_{L^2(\R^d)} \leq C\| S_{\mathcal{A},1} f\|_{H^s(\R^d)} \leq C \| f\|_{H^s(\R^d)}~~~~~{\rm for~all}~~s\in (0,\frac12].
\end{align*}
In particular,  $M_b  S_{\mathcal{A},1}$ is a  lower order operator, and hence, the standard perturbation theory (cf. \cite[Section 2.4]{Lunardi}) implies that $-\mathcal{J}_{\mathcal{A}}$ generates a strongly continuous  and analytic semigroup $\{e^{-t \mathcal{J}_{\mathcal{A}}}\}_{t\geq 0}$ in $L^2 (\R^d)$. It also follows that $\mathcal{J}_{\mathcal{A}}$ admits a bounded $H^\infty$ calculus in $L^2 (\R^d)$. Now we observe that $u(t)=e^{-t\mathcal{J}_{\mathcal{A}}}f$, $f\in L^2 (\R^d)$, satisfies $0=\partial_t u + \mathcal{J}_{\mathcal{A}} u =\partial_t u + \Lambda_{\mathcal{A}} u$ for $t>0$ and $u(t)\rightarrow f$ as $t\rightarrow 0$ in $L^2 (\R^d)$. Thus we also have the representation $u(t) = e^{-t\Lambda_{\mathcal{A}}} f$, that is, $e^{-t\mathcal{J}_{\mathcal{A}}} = e^{-t \Lambda_{\mathcal{A}}}$ for all $t\geq 0$. This proves $\mathcal{J}_{\mathcal{A}}=\Lambda_{\mathcal{A}}$, i.e., $D_{L^2}(\Lambda_{\mathcal{A}}) = D_{L^2}(\mathcal{J}_{\mathcal{A}})=H^1 (\R^d)$. The proof is complete.

\subsection{Domain of Poisson operator in $H^1 (\R^d)$}\label{subsec.domain.H^1}

In this section we study the Poisson operator in $H^1 (\R^d)$. First we consider the operator $\mathcal{Q}_{\mathcal{A}}= M_b (M_{1/\bar{b}} \mathcal{P}_{\mathcal{A}^*})^*$ in $L^2 (\R^d)$, which is the key step to prove $D_{H^1}(\mathcal{P}_{\mathcal{A}}) = H^2 (\R^d)$.

\begin{thm}\label{thm.domain.Q}   Let $\mathcal{Q}_{\mathcal{A}} = M_{1/b} ( M_{\bar{b}} \mathcal{P}_{\mathcal{A}^*})^*$. Then $D_{L^2}(\mathcal{Q}_{\mathcal{A}}) = H^1 (\R^d)$ with equivalent norms and $\mathcal{Q}_{\mathcal{A}} f = M_{1/b} \Lambda_{\mathcal{A}} f - M_{{\bf r_1}/b} \cdot \nabla_x f - M_{(\nabla_x\cdot {\bf r_1})/b} f$ ~ for $f\in H^1 (\R^d)$. 
\end{thm}

\noindent {\it Proof.} Assume that $f\in H^1 (\R^d) = D_{L^2}(\Lambda_{\mathcal{A}})$ (by Theorem \ref{thm.domain.DN}). Then for any $g\in D_{L^2} (\mathcal{P}_{\mathcal{A}^*})= H^1 (\R^d)$ (by Corollary \ref{cor.lem.U_1'}) we have from Proposition \ref{prop.pre.3}, 
\begin{align*}
\langle f, M_{\bar{b}} \mathcal{P}_{\mathcal{A}^*} g\rangle _{L^2 (\R^d)} & = \langle f, \Lambda_{\mathcal{A}^*} g + M_{ {\bf \bar{r}_1}}\cdot \nabla_x g\rangle _{L^2(\R^d)} = \langle \Lambda_{\mathcal{A}} f - \nabla_x \cdot M_{{\bf r_1}} f, g\rangle_{L^2 (\R^d)}.
\end{align*}
In particular, we have the estimate $|\langle f, M_{\bar{b}} \mathcal{P}_{\mathcal{A}^*} g\rangle _{L^2 (\R^d)}|\leq C \| f \|_{H^1 (\R^d)} \| g\|_{L^2(\R^d)}$ for all $g\in D_{L^2} (\mathcal{P}_{\mathcal{A}^*})$. This implies $f\in D_{L^2}( (M_{\bar{b}} \mathcal{P}_{\mathcal{A}^*} )^* )=D_{L^2}( \mathcal{Q}_{\mathcal{A}} )$ and 
\[
\mathcal{Q}_{\mathcal{A}} f =  M_{1/b} \Lambda_{\mathcal{A}} f - M_{{\bf r_1}/b} \cdot \nabla_x f - M_{(\nabla_x\cdot {\bf r_1})/b} f, ~~~~~~~~ f\in H^1 (\R^d),
\]
that is, we have from \eqref{proof.thm.domain.DN.1}-\eqref{proof.thm.domain.DN.2},
\begin{align}
\mathcal{Q}_{\mathcal{A}} f =  -i \mu_{\mathcal{A}} (\cdot,D_x) f  - M_{({\bf r_1} + {\bf r_2})/b} \cdot \nabla_x f - S_{\mathcal{A},1} f - M_{\nabla_x\cdot {\bf r_1}} f, ~~~~~~~~ f\in H^1 (\R^d). \label{proof.thm.domain.Q.1}
\end{align}
To prove the converse embedding we appeal to the argument of the proof of Theorem  \ref{thm.domain.DN}. Set $q_{\mathcal{A}}(x,\xi) =\mu_{\mathcal{A}} (x,\xi) + b(x)^{-1}({\bf r_1} (x) + {\bf r_2} (x) ) \cdot \xi$ and let $\mathcal{K}_{\mathcal{A}}$ be the linear operator in $L^2 (\R^d)$ defined by 
\begin{align}
D_{L^2}(\mathcal{K}_{\mathcal{A}}) = H^1 (\R^d ), ~~~~~~~\mathcal{K}_{\mathcal{A}} f = - i q_{\mathcal{A}} (\cdot,D_x) f  - S_{\mathcal{A},1} f - M_{(\nabla_x\cdot {\bf r_1})/b} f.\label{proof.thm.domain.Q.2}
\end{align}
Then $\mathcal{K}_{\mathcal{A}}f = \mathcal{Q}_{\mathcal{A}}f$ for $f\in H^1 (\R^d)$ by \eqref{proof.thm.domain.Q.1}. On the other hand, Lemma \ref{lem.invariant.lambda} with Remark \ref{rem.pseudo} shows that $ i q_{\mathcal{A}} (\cdot,D_x)$ generates a strongly continuous and analytic semigroup in $L^2 (\R^d)$, and hence, so is true for $\mathcal{K}_{\mathcal{A}}$, since the operators $S_{\mathcal{A},1}$ and  $M_{(\nabla_x\cdot {\bf r_1})/b}$ are of  lower order. Then, as arguing as in the proof of Theorem \ref{thm.domain.DN}, we conclude that $e^{-t\mathcal{K}_{\mathcal{A}}} = e^{-t\mathcal{Q}_{\mathcal{A}}}$ for all $t\geq 0$. Thus we have $\mathcal{K}_{\mathcal{A}} = \mathcal{Q}_{\mathcal{A}}$, as desired. The proof is complete.

\vspace{0.5cm}

\begin{thm}\label{thm.domain.poisson.H^1}  The restriction of  $\{e^{-t\mathcal{P}_{\mathcal{A}}}\}_{t\geq 0}$ in $L^2 (\R^d)$ on the invariant  subspace $H^1 (\R^d)$ defines a strongly continuous  and analytic semigroup. Moreover,  we have $D_{H^1} (\mathcal{P}_{\mathcal{A}}) =  H^2 (\R^d)$ with equivalent norms.  

\end{thm}

\noindent {\it Proof.} The first statement is trivial since we have already proved $D_{L^2}(\mathcal{P}_{\mathcal{A}}) =H^1 (\R^d)$. Thus it suffices to show $D_{H^1}(\mathcal{P}_{\mathcal{A}}) = H^2 (\R^d)$. Since $D_{L^2}(\mathcal{P}_{\mathcal{A}}) = D_{L^2}(\mathcal{P}_{\mathcal{A}^*}) =H^1 (\R^d)$ by Corollary \ref{cor.lem.U_1'}, we have from \eqref{eq.prop.pre.3.3'},
\begin{align*}
u\in H^2 (\R^d) \Longleftrightarrow u\in D_{L^2} (\mathcal{A}') \Longleftrightarrow u\in D_{L^2}(\mathcal{Q}_{\mathcal{A}}\mathcal{P}_{\mathcal{A}} )\Longleftrightarrow \mathcal{P}_{\mathcal{A}} u\in D_{L^2} (\mathcal{Q}_{\mathcal{A}}) \Longleftrightarrow \mathcal{P}_{\mathcal{A}} u\in H^1 (\R^d),
\end{align*}
where we have used $D_{L^2}(\mathcal{Q}_{\mathcal{A}}) = H^1 (\R^d)$ by Theorem \ref{thm.domain.Q}. It is also easy to see that $\|u\|_{H^2(\R^d)}\simeq \| \mathcal{P}_{\mathcal{A}} u \|_{H^1 (\R^d)} + \| u\|_{H^1(\R^d)}$. The proof is complete.

\subsection{Further estimates for remainder part of Poisson operator in $H^s(\R^d)$}\label{subsec.S_A}

Let $S_{\mathcal{A},1}$ be the bounded linear operator in $H^{1/2}(\R^d)$ defined by $S_{\mathcal{A},1}h =\displaystyle \lim_{t\rightarrow 0} \dd/\dd t~U_{\mathcal{A},1} (t) h$  as in the proof of Lemma \ref{lem.U_1'}. In this section we study the mapping property of $S_{\mathcal{A},1}$ in $H^s(\R^d)$.
\begin{prop}\label{prop.S_A} Let $0<s,\epsilon <1$. Then we have 
\begin{align}
\|  S_{\mathcal{A},1} h \|_{H^s(\R^d)}\leq C \| h \|_{H^s(\R^d)},~~~~~~~~\|  S_{\mathcal{A},1} h \|_{H^1(\R^d)}\leq C \| h \|_{H^{1+\epsilon}(\R^d)}. \label{est.prop.S_A.1}
\end{align}
\end{prop}

\noindent {\it Proof.} Let $h\in \mathcal{S}(\R^d)$.  The estimate \eqref{est.prop.S_A.1} with $s\in (0,1/2]$ is already proved by \eqref{proof.lem.U_1'.5}. Next we consider the case $s\in (1/2,1)$. Let us recall that $U_{\mathcal{A},1}(t)h$ is the solution to \eqref{eq.dirichlet} with $F$ given by \eqref{proof.lem.U_1.1} and $g=0$. By \cite[Theorem 5.1]{MaekawaMiura1} the characterization of $D_{L^2}(\mathcal{P}_{\mathcal{A}}) = D_{L^2}(\mathcal{P}_{\mathcal{A}^*}) = H^1(\R^d)$ provides the integral representation of $U_{\mathcal{A},1}(t)h$ such that 
\begin{align*}
U_{\mathcal{A},1} (t) h & = \int_0^t e^{-(t-s)\mathcal{P}_{\mathcal{A}}} \int_s^\infty e^{-(\tau-s)\mathcal{Q}_{\mathcal{A}}} M_{1/b}\big (  \nabla\cdot M_\chi \Pi h + M_\chi G_\zeta h + R h\big ) \dd\tau \dd s,
\end{align*}
which gives 
\begin{align}
S_{\mathcal{A},1} h = \int_0^\infty e^{-\tau \mathcal{Q}_{\mathcal{A}}}  M_{1/b}\big ( \nabla\cdot M_\chi \Pi h + M_\chi G_\zeta h + R h\big ) \dd \tau. \label{proof.prop.S_A.1}
\end{align}
Here we will only show 
\begin{align}
\| \int_0^\infty  e^{-\tau \mathcal{Q}_{\mathcal{A}}}  M_{1/b}  \nabla_x \cdot M_\chi \Pi' h \dd \tau \|_{\dot{H}^s(\R^d)} \leq C \| h \|_{H^s(\R^d)}, \label{proof.prop.S_A.2}
\end{align}
for the other terms are treated in the similar manner. We note that the term $\Pi'h$ is not differentiable in $x$ (see the definition \eqref{proof.lem.U_1.0}), and thus, the term $e^{-\tau \mathcal{Q}_{\mathcal{A}}} M_{1/b}  \nabla_x \cdot M_\chi \Pi' h$ in \eqref{proof.prop.S_A.2} has to be interpreted as  
\begin{align}
e^{-\tau \mathcal{Q}_{\mathcal{A}}} M_{1/b}  \nabla_x  \cdot M_\chi \Pi' h & = ( I + \mathcal{Q}_{\mathcal{A}} ) e^{-\tau\mathcal{Q}_{\mathcal{A}}} \big (- \nabla_x ( I + \mathcal{P}_{\mathcal{A}^*} )^{-1} M_{1/\bar{b}}\big )^* \cdot M_\chi \Pi'h,\label{proof.prop.S_A.3}
\end{align}
where we have used the formal adjoint relation $( I + \mathcal{Q}_{\mathcal{A}} )^{-1} M_{1/b}  \nabla_x = \big (- \nabla_x ( I + \mathcal{P}_{\mathcal{A}^*} )^{-1} M_{1/\bar{b}}\big )^*$. Since $\big (- \nabla_x ( I + \mathcal{P}_{\mathcal{A}^*} )^{-1} M_{1/\bar{b}}\big )^*$ is a bounded linear operator in $L^2 (\R^d)$ by $D_{L^2}(\mathcal{P}_{\mathcal{A}^*})=H^1 (\R^d)$, the right-hand side of \eqref{proof.prop.S_A.3} is well-defined for each $\tau>0$. Then from $\mathcal{Q}_{\mathcal{A}}e^{-\tau \mathcal{Q}_{\mathcal{A}}} = - \dd/\dd\tau e^{-\tau \mathcal{Q}_{\mathcal{A}}}$ and from the integration by parts together with $\Pi 'h|_{t=0} =0$ for $h\in \mathcal{S}(\R^d)$ (due to the definition of $\Pi'$) the estimate \eqref{proof.prop.S_A.2} is essentially reduced to 
\begin{align}
\| \int_0^\infty  e^{-\tau \mathcal{Q}_{\mathcal{A}}}  B^* \cdot M_\chi \partial_\tau \Pi'h\dd \tau \|_{\dot{H}^s(\R^d)} \leq C \| h \|_{H^s(\R^d)},~~~ ~~~ B=- \nabla_x ( I + \mathcal{P}_{\mathcal{A}^*} )^{-1} M_{1/\bar{b}},\label{proof.prop.S_A.4}
\end{align}
since the other terms are of lower order. We appeal to the duality argument and consider the integral 
\begin{align}
& ~~~\langle (-\Delta_x)^\frac{s}{2}\int_0^\infty  e^{-\tau \mathcal{Q}_{\mathcal{A}}}  B^* \cdot M_\chi \partial_\tau \Pi'h\dd \tau,  \varphi \rangle _{L^2 (\R^d)} \nonumber \\
& = \int_0^\infty \langle  B^* \cdot M_\chi \partial_\tau \Pi'h, M_{\bar{b}} e^{-\tau \mathcal{P}_{\mathcal{A}^*}} M_{1/\bar{b}} (-\Delta_x)^\frac{s}{2} \varphi \rangle _{L^2 (\R^d)} \dd\tau\label{proof.prop.S_A.5}
\end{align}
for $\varphi\in \mathcal{S}(\R^d)$. Then  we have 
\begin{align}
&~~~ {\rm R.H.S.~of~\eqref{proof.prop.S_A.5}}\nonumber \\
&  \leq C \big ( \int_0^\infty M_\chi \tau^{2(1-s)}\| \partial_\tau \Pi'h \|_{L^2 (\R^d)}^2 \frac{\dd \tau}{\tau} \big )^\frac12 \big ( \int_0^\infty M_\chi \tau^{2 s}\|  M_{\bar{b}} e^{-\tau \mathcal{P}_{\mathcal{A}^*}} M_{1/\bar{b}} (-\Delta_x)^\frac{s}{2} \varphi \|_{L^2 (\R^d)}^2\frac{\dd \tau}{\tau} \big )^\frac12.\label{proof.prop.S_A.6}
\end{align}
By \eqref{proof.lem.U_1.0}  we see $\partial_\tau \Pi'= G_{i (1 + \tau i\mu_{\mathcal{A}})A'\nabla_x\mu_{\mathcal{A}} }$ and it is straightforward to check that 
\[
p (x,\xi,\tau) = \tau^{1-s} |\xi|^{-s}(1 + \tau i\mu_{\mathcal{A}})  i A'\nabla_x\mu_{\mathcal{A}}, ~~~~~~~s\in (0,1)
\]
satisfies the condition \eqref{assume.lem.G_p.2}. Hence we have from \eqref{est.lem.G_p.2},
\begin{align*}
\int_0^\infty M_\chi \tau^{2(1-s)}\| \partial_\tau \Pi' h \|_{L^2 (\R^d)}^2 \frac{\dd \tau}{\tau}  = \int_0^\infty M_\chi \| G_p (-\Delta_x)^\frac{s}{2} h \|_{L^2 (\R^d)}^2 \frac{\dd\tau}{\tau}\leq C \| (-\Delta_x)^\frac{s}{2} h \|_{L^2 (\R^d)}^2.
\end{align*}
Next we estimate the second integral of the right-hand side of \eqref{proof.prop.S_A.6}. By the duality argument and $D_{L^2}(\mathcal{Q}_{\mathcal{A}}) = H^1 (\R^d)$ it is easy to see $\| M_{\bar{b}} e^{-\tau \mathcal{P}_{\mathcal{A}^*}} M_{1/\bar{b}} (-\Delta_x)^\frac{\kappa}{2} \varphi \|_{L^2 (\R^d)}\leq C \tau^{-\kappa} \|\varphi \|_{L^2 (\R^d)}$ for any $\kappa\in [0,1]$ and $\tau\in (0,2)$. Let $\{\psi_r\}_{r>0}$ be the family of functions introduced in Appendix \ref{appendix.proof.lemma} (with $s$ replaced by $r$). Then one can verify the estimates
\begin{align*}
\tau^s \| M_{\bar{b}} e^{-\tau \mathcal{P}_{\mathcal{A}^*}} M_{1/\bar{b}} (-\Delta_x)^\frac{s}{2}\psi_r * \varphi \|_{L^2 (\R^d)} & =\tau^s \| M_{\bar{b}} e^{-\tau \mathcal{P}_{\mathcal{A}^*}} M_{1/\bar{b}} (-\Delta_x)^\frac{1}{2} (-\Delta_x)^\frac{s-1}{2}\psi_r * \varphi \|_{L^2 (\R^d)}\\
& \leq C \tau^{-1+s}r^{1-s} \|\varphi\|_{L^2(\R^d)}, ~~~~~~~~~~~0<r\leq \tau<2,\\
\tau^s \| M_{\bar{b}} e^{-\tau \mathcal{P}_{\mathcal{A}^*}} M_{1/\bar{b}} (-\Delta_x)^\frac{s}{2}\psi_r * \varphi \|_{L^2 (\R^d)} & \leq \tau^s \|(-\Delta_x)^\frac{s}{2}\psi_r * \varphi \|_{L^2 (\R^d)}\\
& \leq C \tau^{s}r^{-s} \|\varphi\|_{L^2(\R^d)}, ~~~~~~~~~~~0<\tau\leq r, ~~ 0<\tau<2.
\end{align*}
Hence the Schur lemma \cite[pp.643-644]{Grafakos} yields
\begin{align*}
\int_0^\infty M_\chi \tau^{2 s}\|  M_{\bar{b}} e^{-\tau \mathcal{P}_{\mathcal{A}^*}} M_{1/\bar{b}} (-\Delta_x)^\frac{s}{2} \varphi \|_{L^2 (\R^d)}^2\frac{\dd \tau}{\tau} \leq C \| \varphi \|_{L^2 (\R^d)}^2,
\end{align*}
as desired. This completes the proof of \eqref{est.prop.S_A.1} with $s\in (0,1)$. To prove the second estimate in \eqref{est.prop.S_A.1}  we go back to the representation \eqref{proof.prop.S_A.1}. Here we will only show, instead of \eqref{proof.prop.S_A.4},
\begin{align}
\| \int_0^\infty  e^{-\tau \mathcal{Q}_{\mathcal{A}}}  B^* \cdot M_\chi \partial_\tau \Pi' h\dd \tau \|_{\dot{H}^1(\R^d)} \leq C \| h \|_{H^{1+\epsilon} (\R^d)}, ~~~~~~\epsilon\in (0,1).\label{proof.prop.S_A.7}
\end{align}
We use the identity $\partial_\tau \Pi' h = G_q (I-\Delta_x)^{(1+\epsilon)/2} h$ with $q= i (1+|\xi|^2)^{-(1+\epsilon)/2}
 (1+ \tau i\mu_{\mathcal{A}})A'\nabla_x \mu_{\mathcal{A}} $. For $h\in \mathcal{S}(\R^d)$ the limit $\displaystyle \lim_{\tau \rightarrow 0}G_q (I-\Delta_x)^{(1+\epsilon)/2} h$ exists in $L^2 (\R^d)$ and we also have 
\begin{align}
\sup_{0<\tau<2} \| G_q (\tau)  f \|_{L^2 (\R^d)}\leq C \| f\|_{L^2 (\R^d)},~~~~~~~f\in \mathcal{S}(\R^d ). \label{proof.prop.S_A.8}
\end{align}
The proof of \eqref{proof.prop.S_A.8} is postponed to the appendix. Then from $D_{L^2}(\mathcal{Q}_{\mathcal{A}}) = H^1 (\R^d)$ we have 
\begin{align}
{\rm L.H.S.~of~\eqref{proof.prop.S_A.7}} &\leq C \| \mathcal{Q}_{\mathcal{A}}\int_0^\infty M_\chi e^{-\tau \mathcal{Q}_{\mathcal{A}}} B^* \cdot G_q (I-\Delta_x)^{\frac{1+\epsilon}{2}} h \dd\tau \|_{L^2 (\R^d)} + ~{\rm (lower~order)} \nonumber \\
& \leq C \lim_{\tau\rightarrow 0} \| B^* \cdot G_q (I-\Delta_x)^{\frac{1+\epsilon}{2}} h \|_{L^2 (\R^d)} \nonumber  \\
& ~~~~~~~~~~ + C \int_0^\infty M_\chi \| e^{-\tau \mathcal{Q}_{\mathcal{A}}} B^* \cdot \partial_\tau G_q (I-\Delta_x)^{\frac{1+\epsilon}{2}} h  \|_{L^2 (\R^d)} \dd\tau +~{\rm (lower ~order)}\nonumber \\
&\leq C \| h \|_{H^{1+\epsilon}(\R^d)} + C \int_0^\infty M_\chi \| \partial_\tau G_q (I-\Delta_x)^{\frac{1+\epsilon}{2}} h \|_{L^2 (\R^d)} \dd\tau +{\rm ~ (lower~order)}.  \label{proof.prop.S_A.9}
\end{align}
By the definition of $q$ we see $\partial_\tau G_q = \tau^{-1+\epsilon} G_{\tilde q}$ with $\tilde q =\tau^{1-\epsilon} ( i q\mu_{\mathcal{A}}  + \partial_\tau q )$ and $\tilde q$ satisfies the condition \eqref{assume.lem.G_p.1} with $T=2$. Thus, \eqref{est.lem.G_p.1} implies ~${\rm R.H.S.~of~\eqref{proof.prop.S_A.9}}\leq C \| h \|_{H^{1+\epsilon}(\R^d)}$. The proof is complete.

\subsection{Proof of Theorems \ref{thm.main.1}, \ref{thm.main.2}, and \ref{thm.main.3}}\label{subsec.proof.thm}

\noindent {\it Proof of Theorem \ref{thm.main.1}.} The assertion (ii) of Theorem \ref{thm.main.1} with $s=0$ and $s=1$ is already proved in Corollary \ref{cor.lem.U_1'} and Theorem \ref{thm.domain.poisson.H^1}. Then the case $s\in (0,1)$ follows from the interpolation inequality and the details are omitted. It remains to show the last statement of (i). By Theorem \ref{thm.domain.poisson.H^1} we have 
\begin{align}
u\in D_{L^2}(\mathcal{P}_{\mathcal{A}}^2) \Longleftrightarrow  \mathcal{P}_{\mathcal{A}}u \in D_{L^2} (\mathcal{P}_{\mathcal{A}}) = H^1 (\R^d) \Longleftrightarrow  u\in D_{H^1} (\mathcal{P}_{\mathcal{A}}) = H^2 (\R^d).
\end{align}
It is also easy to see the norm equivalence between $H^2 (\R^d)$ and $D_{L^2} (\mathcal{P}_{\mathcal{A}}^2 )$. Then the sectorial operator $T=I + \mathcal{P}_{\mathcal{A}}$ in $L^2 (\R^d)$, which is invertible by \cite[Remark 2.6]{MaekawaMiura1},  satisfies 
\begin{align}
(L^2 (\R^d), D_{L^2} (T^2) )_{\frac12,2} = (L^2 (\R^d), H^2 (\R^d) )_{\frac12,2} = H^1 (\R^d) = D_{L^2} (T).\label{identity.komatsu}
\end{align}
By the Komatsu theorem \cite[Theorem 6.6.8]{Haase} the identity \eqref{identity.komatsu} implies that the operator $T$ admits a bounded $H^\infty$ calculus in $L^2(\R^d)$. The proof is complete.

\vspace{0.5cm}

\noindent {\it Proof of Theorem \ref{thm.main.2}.} The assertions follow from Corollary \ref{cor.lem.U_1'} and Theorem \ref{thm.domain.Q} together with Proposition \ref{prop.pre.3}. The proof is complete.

\vspace{0.5cm}

\noindent {\it Proof of Theorem \ref{thm.main.3}.} The case $s=0$ is already proved by Theorem \ref{thm.domain.DN}. It suffices to consider the endpoint case $s=1$.  Theorem \ref{thm.domain.DN} implies that $H^1 (\R^d)$ is invariant under the action of $\{e^{-t\Lambda_{\mathcal{A}}}\}_{t\geq 0}$ and the restriction of this semigroup in $H^1 (\R^d)$ is also analytic and strongly continuous. Hence it suffices to show that the generator of this restriction semigroup satisfies $D_{H^1}(\Lambda_{\mathcal{A}}) = H^2(\R^d)$ with equivalent norms. By the proof of Theorem \ref{thm.domain.DN} we have 
\begin{align*}
\Lambda_{\mathcal{A}}=\mathcal{J}_{\mathcal{A}},~~~~~~~\mathcal{J}_{\mathcal{A}}=-i\lambda_{\mathcal{A}}(\cdot,D_x) - M_b S_{\mathcal{A},1}~~~~~{\rm as ~a ~operator ~in~}~ L^2 (\R^d),
\end{align*}
where $\lambda_{\mathcal{A}}(\cdot,D_x)$ is the pseudo-differential operator with its symbol $\lambda_{\mathcal{A}}(x,\xi) = b(x)\mu_{\mathcal{A}}(x,\xi) +{\bf r_2} (x) \cdot \xi$. On the other hand, by Lemma \ref{lem.invariant.lambda} and Theorem \ref{thm.class.Phi} the operator $i\lambda_{\mathcal{A}}(\cdot,D_x)$ generate a strongly continuous and analytic semigroup in $H^1(\R^d)$, and $D_{H^1}(\lambda_{\mathcal{A}}(\cdot,D_x)) = H^2(\R^d)$ holds with equivalent norms. Then so is true for $\mathcal{J}_{\mathcal{A}}$, for $b$ is Lipschitz and $S_{\mathcal{A},1}$ is of lower order by Proposition \ref{prop.S_A}. Since it is easy to see that $\Lambda_{\mathcal{A}}=\mathcal{J}_{\mathcal{A}}$ as a operator in $H^1 (\R^d)$, we conclude that $D_{H^1}(\Lambda_{\mathcal{A}}) = H^2(\R^d)$ holds with equivalent norms. The proof is complete.

\begin{rem}
\label{Hinfty}
{\rm In the proof of Theorem \ref{thm.main.1} we have  established the expansion $\mathcal{P}_{\mathcal{A}}=-i \mu_{\mathcal{A}} (\cdot, D_x) + R$, where $\mu_{\mathcal{A}}(\cdot, D_x)$ is the pseudo-differential operator with symbol \eqref{mu_A.intro} and $R$ is a bounded operator in $H^s(\R^d)$ for $s\in (0,1)$, while it is a bounded operator from $H^{1+\epsilon}(\R^d)$, $\epsilon>0$, to $H^1 (\R^d)$. Similar expansion is obtained also for $\Lambda_{\mathcal{A}}$. Then one can apply the results of \cite[Theorem 4.8]{ES} to obtain a stronger statement that $\mathcal{P}_{\mathcal{A}}$ (and $\Lambda_{\mathcal{A}}$)  admits a  bounded $H^\infty$ calculus in $H^s(\R^d)$, $s\in [0,1)$. 
Our proof of bounded $H^\infty$ calculus for $\mathcal{P}_{\mathcal{A}}$ is based on the characterization $D_{H^1}(\mathcal{P}_{\mathcal{A}})=H^2 (\R^d)$ and the Komatsu theorem \cite[Proposition 2.7]{Haase},  which is different from the approach in \cite{ES}.
}
\end{rem}
\appendix 
\section{Appendix}\label{sec.appendix}

\subsection{Remark on  pseudo-differential operator $\mu_{\mathcal{A}}(\cdot,D_x)$}\label{appendix.pseudo-differential}

In view of the definition  \eqref{def.mu_A} for $\mu_{\mathcal{A}}(x,\xi)$, which is the root of \eqref{eq.root} with positive imaginary part, it is natural to introduce the map $\Phi: A\longmapsto \mu_{\mathcal{A}}$, i.e.,
\begin{align}
\Phi (A) =- \frac{{\bf v} (x) \cdot \xi}{2} + i \big \{ \frac{1}{b (x) } \langle A'(x) \xi,\xi\rangle -\frac{1}{4}({\bf v}(x) \cdot \xi )^2\big \}^\frac12, ~~~~~~~~~{\bf v} = \frac{{\bf r_1} + {\bf r_2}}{b}, \label{def.appendix.Phi}
\end{align}
where $A$ is a matrix satisfying the ellipticity condition \eqref{ellipticity}. We denote by $\mathcal{R}_{Lip} (\Phi)$ the range of $\Phi$ for the Lipschitz class of $A$, that is,
\begin{align*}
\mathcal{R}_{Lip} (\Phi) & = \{ \mu (x,\xi)\in {\rm Lip} (\R^d \times \R^d)~|~{\rm there~ is ~a ~matrix}~A\in ({\rm Lip}(\R^d ) )^{(d+1)\times (d+1)}~{\rm satisfying }~\eqref{ellipticity}\nonumber \\
& ~~~~~~~~~~~~~~~~~~~~~~~~~~~~~~~~~~~~~~~ {\rm for ~ some~} \nu_1,\nu_2>0~{\rm such ~that~} \mu =\Phi (A)\}. 
\end{align*}

\noindent The next lemma is used in the study of $\Lambda_{\mathcal{A}}$ and $\mathcal{Q}_{\mathcal{A}}$. 
\begin{lem}\label{lem.invariant.lambda} Assume that $\mu = \Phi (A) \in \mathcal{R}_{Lip} (\Phi)$ with $A=(a_{i,j})_{1\leq i,j\leq d+1}$.  Set  $b=a_{d+1,d+1}$, ${\bf r_1} =  (a_{j,d+1})_{1\leq j\leq d}$, and ${\bf r_2} = (a_{d+1,j})_{1\leq j\leq d}$.  Then the functions $\lambda(x,\xi)$ and $q (x,\xi)$ defined by 
\begin{align}
\lambda (x,\xi ) = b (x) \mu (x,\xi) + {\bf r_2} (x) \cdot \xi,  ~~~~~~~ q (x,\xi) = \mu (x,\xi) + \frac{{\bf r_1} (x) + {\bf r_2} (x) }{b(x)} \cdot \xi,
\end{align}
belong to $ \mathcal{R}_{Lip} (\Phi)$.

\end{lem}

\noindent {\it Proof.} Since $\mu$ solves \eqref{eq.root}, $\lambda$ and $q$ respectively satisfy
\begin{align*}
&\lambda^2 + ({\bf r_1} (x) - {\bf r_2} (x) )\cdot \xi \lambda + b(x) \langle A' (x) \xi,\xi\rangle - {\bf r_1}(x) \cdot \xi ~{\bf r_2}(x) \cdot \xi = 0,\\
& b(x) q^2 - ({\bf r_1} (x) + {\bf r_2} (x) )\cdot \xi q + \langle A' (x) \xi, \xi \rangle = 0.
\end{align*}
Set $M' = \big (m_{i,j} \big )_{1\leq i,j\leq d} = \big (|b|^2 a_{i,j} - \bar{b}a_{d+1,i} a_{j,d+1} \big )_{1\leq i,j\leq d}$, ${\bf s_1} = \bar{b}{\bf r_1}$, ${\bf s_2}= -\bar{b}{\bf r_2}$, and set $N'=A'$, ${\bf u_1}=-{\bf r_1}$, ${\bf u_2}=-{\bf r_2}$. Then  the matrices 
\begin{equation}
M  =
\begin{pmatrix}
\mbox{} & \mbox{} & \mbox{} & \mbox{} \\
\mbox{} & \large{M'} & \mbox{} &  {\bf s_1} \\
\mbox{} & \mbox{} & \mbox{} & \mbox{}\\
\mbox{} & {\bf s_2}^\top & \mbox{} & \bar{b} \end{pmatrix}, ~~~~~~~~~N  =
\begin{pmatrix}
\mbox{} & \mbox{} & \mbox{} & \mbox{} \\
\mbox{} & \large{N'} & \mbox{} &  {\bf u_1} \\
\mbox{} & \mbox{} & \mbox{} & \mbox{}\\
\mbox{} & {\bf u_2}^\top & \mbox{} & b \end{pmatrix}~
\end{equation}
satisfy \eqref{ellipticity} and \eqref{lipschitz} (with possibly different ellipticity constants), while we have $\lambda=\Phi (M)$ and $q=\Phi (N)$. The proof is complete.

\vspace{0.5cm}

As stated in Remark \ref{rem.pseudo}, in the proof of the characterization $D_{H^s}(\mathcal{P}_{\mathcal{A}}) = H^{1+s} (\R^d)$ we did not use the mapping property of  the pseudo-differential operator $\mu_{\mathcal{A}} (\cdot, D_x)$, which has a representation 
\begin{align}
\mu_{\mathcal{A}} (\cdot,  D_x ) f (x) =  \frac{1}{(2\pi)^\frac{d}{2}} \int_{\R^d} \mu_{\mathcal{A}} (x,\xi )  e^{ix\cdot \xi} \hat{f} (\xi) \dd \xi ~~~~~~~~{\rm for}~~f\in \mathcal{S} (\R^d).
\end{align}

\noindent The results of Theorem \ref{thm.main.1} and Proposition \ref{prop.S_A} yield 
\begin{thm}\label{thm.class.Phi} Let $s\in [0,1]$. Let $\mu = \Phi (A) \in \mathcal{R}_{Lip} (\Phi)$.  Then the associated pseudo-differential operator  $i\mu (\cdot,D_x)$ generate a strongly continuous and analytic semigroup in  $H^s(\R^d)$, and $D_{H^s}(\mu (\cdot, D_x) ) = H^{1+s}(\R^d)$ holds with equivalent norms. Moreover, for the Poisson operator $\mathcal{P}_{\mathcal{A}}$ associated with $\mathcal{A}=-\nabla\cdot A\nabla$ we have the identity
\begin{align}
\mathcal{P}_{\mathcal{A}} f = - i \mu (\cdot, D_x) f- S_{\mathcal{A},1}  f,~~~~~~~~f\in H^{1} (\R^d),\label{eq.thm.class.Phi}
\end{align}
where $S_{\mathcal{A},1}$ is the linear operator defined in the proof of Lemma \ref{lem.U_1'}, which is bounded in $H^s(\R^d)$, $s\in (0,1)$. In particular, $-i\mu (\cdot,D_x)$ admits a bounded $H^\infty$ calculus in $L^2 (\R^d)$.

\end{thm}

\begin{rem}{\rm In fact, by applying the general results of \cite{ES} for pseudo-differential operators with nonsmooth symbols, it follows that $-i\mu (\cdot,D_x)$ admits a bounded $H^\infty$ calculus in $H^s (\R^d)$, $s\in [0,1)$. 
In this sense,  the properties of $i\mu (\cdot,D_x)$ stated in Theorem \ref{thm.class.Phi} themselves  are not essentially new. As commented in Remark \ref{rem.pseudo}, the special feature of our proof is that we use the information of $\mathcal{P}_{\mathcal{A}}$ to derive the properties of  $-i\mu (\cdot,D_x)$, where the underlying key structure is the factorizations of $\mathcal{A}'$ and $\mathcal{A}$ in Theorem \ref{thm.main.2}. 
 
}
\end{rem}
 
\noindent {\it Proof of Theorem \ref{thm.class.Phi}.} For $f\in H^{1+s}(\R^d)$ we define $\mu (\cdot, D_x) f =  \displaystyle \lim_{n\rightarrow \infty}\mu (\cdot,D_x)f_n$ in $H^s(\R^d)$, where $\{f_n\}$ is a sequence in $\mathcal{S}(\R^d)$ converging to $f$ in $H^{1+s}(\R^d)$. This is well-defined since  \eqref{eq.thm.class.Phi} holds for $f\in \mathcal{S}(\R^d)$, and then Theorem \ref{thm.main.1} and Proposition \ref{prop.S_A} imply $\|\mu (\cdot, D_x) f \|_{H^s(\R^d)} \leq \| \mathcal{P}_{\mathcal{A}} f \|_{H^s(\R^d)} + \| S_{\mathcal{A},1} f \|_{H^s(\R^d)} \leq C \| f \|_{H^{1+s}(\R^d)}$ for $f\in \mathcal{S}(\R^d)$. Since $D_{H^s}(\mathcal{P}_{\mathcal{A}})=H^{1+s}(\R^d)$ and $\mathcal{P}_{\mathcal{A}}$ is closed in $H^s(\R^d)$, we observe also from Proposition \ref{prop.S_A} that the above realization of $i\mu (\cdot, D_x)$ in $H^s(\R^d)$ satisfies \eqref{eq.thm.class.Phi} for any $f\in H^{1+s}(\R^d)$. Hence $i\mu (\cdot, D_x)$ defined above is a perturbation from $-\mathcal{P}_{\mathcal{A}}$ by  $S_{\mathcal{A},1}$ which is a lower order operator, and the desired properties of $i\mu (\cdot, D_x)$ then follow from the ones of  $-\mathcal{P}_{\mathcal{A}}$ by the standard perturbation theory of sectorial operators. The proof is complete.

\subsection{Proofs of Lemma \ref{lem.G_p} and \eqref{proof.prop.S_A.8}}\label{appendix.proof.lemma}

\noindent {\it Proof of Lemma \ref{lem.G_p}.} We set 
\begin{align}
G_p (x,y,t) = \frac{1}{(2 \pi)^\frac{d}{2}} \int_{\R^d} p (x,\xi, t) e^{i t \mu_{\mathcal{A}} (x,\xi,t )} e^{i y\cdot \xi }\dd \xi .\label{proof.lem.G_p.1}
\end{align}
Then $\big (G_p (t) h \big )(x) = (2\pi)^{d/2} \int_{\R^d} G_p (x,x-y,t) h (y) \dd y$, and thus, it suffices to show 
\begin{align}
|G_p (x,y,t )|\leq C t^{-d} ( 1+\frac{|y|}{t})^{-d-\epsilon},~~~~~~~~0<t<T,~x,y\in \R^d,\label{proof.lem.G_p.2}
\end{align}
for some $\epsilon>0$. When $|y|\leq t$ we have from \eqref{estimate.mu_A.1} and \eqref{assume.lem.G_p.1},
\begin{align*}
|G_p (x,y,t) |\leq C \int_{\R^d} e^{-c t|\xi|} \dd\xi\leq C t^{-d}\leq C t^{-d}  ( 1+\frac{|y|}{t})^{-d-1}.
\end{align*}
Next we consider the case $|y|\geq t$. For any multi-index $\alpha$ with its length $|\alpha|=j_0$, the integration by parts yields
\begin{align*}
y^\alpha G_p (x,y,t) & = \frac{i^{j_0}}{(2 \pi)^\frac{d}{2}} \int_{\R^d} ( \partial_\xi^\alpha p ) e^{it \mu_{\mathcal{A}}} e^{i y\cdot \xi} \dd \xi + \frac{i^{j_0}}{(2 \pi)^\frac{d}{2}} \sum_{\alpha\geq \beta,|\beta|\ne 0} C_{\alpha,\beta} \int_{\R^d} (\partial_\xi ^{\alpha-\beta} p ) ( \partial_\xi^\beta e^{it\mu_{\mathcal{A}}}) e^{i y\cdot\xi} \dd \xi \\
& =: I + II.
\end{align*} 
Let $\chi_R (\xi)$ be a smooth function such that $\chi_R =1$ for $|\xi|\leq R$, $\chi_R =0$ for $|\xi|\geq 2 R$, and $\|\nabla^k \chi_R \|_{L^\infty}\leq C R^{-k}$. We divide $I$ into a low frequency part $I_1$ and a high frequency part $I_2$, where 
\begin{align*}
I_1 =   \frac{i^{j_0}}{(2 \pi)^\frac{d}{2}} \int_{\R^d} \chi_R \cdots \dd \xi, ~~~~~~~~ I_2 =   \frac{i^{j_0}}{(2 \pi)^\frac{d}{2}} \int_{\R^d} (1-\chi_R) \cdots  \dd \xi.
\end{align*}
Then the condition \eqref{assume.lem.G_p.1} leads to $|I_1|\leq C t^{l_{j_0}} \int_{|\xi|\leq 2 R} |\xi|^{-j_0+l_{j_0}} \dd\xi \leq C t^{l_{j_0}} R^{d-j_0 +l_{j_0}}$, while the integration by parts combined with \eqref{estimate.mu_A.1} and \eqref{assume.lem.G_p.1} gives $|y^\gamma I_2|\leq C \int_{|\xi|\geq R} |\xi|^{-d-1}\dd\xi \leq C R^{-1}$ for any multi-index $\gamma$ satisfying $|\gamma|=d+1-j_0$. Similarly, we divide $II$ into  $II_1$ and  $II_2$, where 
\begin{align*}
II_1 =   \frac{i^{j_0}}{(2 \pi)^\frac{d}{2}} \int_{\R^d} \chi_R \cdots \dd \xi, ~~~~~~~~ II_2 =   \frac{i^{j_0}}{(2 \pi)^\frac{d}{2}} \int_{\R^d} (1-\chi_R) \cdots  \dd \xi.
\end{align*}
Then we have from $|\beta|\geq 1$ that $|II_1|\leq C t \int_{|\xi|\leq 2 R} |\xi|^{-j_0 + 1} \dd\xi \leq C t  R^{d-j_0 + 1}$, while $II_2$ is estimated as $|y^\gamma II_2|\leq C R^{-1}$. Collecting these, we see
\begin{align*}
|G_p (x,y,t)|\leq C |y|^{-j_0} (I_1  + I_2 + II_1 + II_2) \leq  C |y|^{-j_0}\big ( t^{l_{j_0}} R^{d-j_0 + l_{j_0}} +  t  R^{d-j_0 +1} + |y|^{-d-1+j_0}R^{-1}\big ).
\end{align*}
By taking 
\begin{align*}
R=t^{-\frac{l_{j_0}}{d+1-j_0+l_{j_0}}} |y|^{-\frac{d+1-j_0}{d+1-j_0+l_{j_0}}}~~~{\rm if}~~l_{j_0}\in (0,1],~~~~~~~~ R=t^{-\frac{1}{d+2-j_0}} |y|^{-\frac{d+1-j_0}{d+2-j_0}}~~~{\rm if}~~l_{j_0}\geq 1,
\end{align*}
we get the desired estimate \eqref{proof.lem.G_p.2} also for $|y|\geq t$. Now the proof of \eqref{est.lem.G_p.1} is complete. To prove \eqref{est.lem.G_p.2} let $\psi\in C_0^\infty (\R^d)$ be a real-valued function with zero average such that 
\[
\int_0^\infty \| \psi_s * f \|_{L^2(\R^d)}^2 \frac{\dd s}{s}  = \| f\|_{L^2 (\R^d)}^2,~~~~~~~~ f\in L^2 (\R^d).
\]
Here $\psi_s (x) = s^{-d} \psi (x/s)$. We may take $\psi=\Delta \tilde \psi$ so that $\| s^{-1} \nabla_x \tilde \psi_s * f \|_{L^2(\R^d)}\leq C \| f\|_{L^2 (\R^d)}$ holds. Thanks to \eqref{est.lem.G_p.1} we have $\| G_p (t) \psi_s* h \|_{L^2(\R^d)}\leq C \| h \|_{L^2 (\R^d)}$ for all $t,s>0$. Moreover, when $t\geq s>0$ we apply  \eqref{est.lem.G_p.1} to $p$ replaced by $p_1 =t\xi p$  and get
\begin{align*}
\| G_p (t) \psi_s*h \|_{L^2(\R^d)} & =  \| G_{p} (t) \nabla_x \cdot \nabla_x \tilde \psi_s * h \|_{L^2(\R^d)} \\
&  = t^{-1} \| G_{p_1} (t) \cdot \nabla_x \tilde \psi_s*h \|_{L^2(\R^d)} \leq C t^{-1} \| \nabla_x \tilde \psi_s * h\|_{L^2(\R^d)} \leq C t^{-1} s\| h \|_{L^2 (\R^d)}.
\end{align*}
When $s\geq t>0$ we take $l=\min \{l_0,\cdots,l_{d+1}\}>0$ and set $p_2 = (t|\xi |)^{-l/2} p$. Then it is easy to see that $p_2$ satisfies \eqref{assume.lem.G_p.1}. Hence we have 
\begin{align*}
\|G_p (t) \psi_s*h \|_{L^2(\R^d)} & = t^{\frac{l}{2}} \| G_{p_2} (t) (-\Delta_x)^\frac{l}{4} \psi_s * h \|_{L^2(\R^d)} \\
& \leq C t^{\frac{l}{2}} \|(-\Delta_x)^{\frac{l}{4}} \psi_s * h\|_{L^2(\R^d)} \leq C t^\frac{l}{2} s^{-\frac{l}{2}}\| h \|_{L^2 (\R^d)}.
\end{align*}
Now we can apply  the Schur lemma (cf. see \cite[pp.643-644]{Grafakos}) to $\{G_p(t)\}_{t>0}$, which leads to \eqref{est.lem.G_p.2}. The proof of Lemma \ref{lem.G_p} is complete.
 
\vspace{0.5cm}

\noindent {\it Proof of \eqref{proof.prop.S_A.8}.} With the notation \eqref{proof.lem.G_p.1} it suffices to show
\begin{align}
|G_q(x,y,t)|\leq C \min\{ |y|^{-d+\delta},~|y|^{-d-\delta}\},~~~~~~~x,y\in \R^d,~0<t<2\label{proof.lem.G_p.3}
\end{align}  
for some $\delta>0$, where $q(x,\xi,t) = i (1+|\xi|^2)^{-(1+\epsilon)/2} (1+i t \mu_{\mathcal{A}} )A'\nabla_x\mu_{\mathcal{A}}$. For any multi-index $\alpha$ with $d-1\leq |\alpha|\leq d$ we have 
\begin{align*}
y^\alpha G_q (x,y,t) & = \frac{i^{|\alpha|}}{(2\pi)^\frac{d}{2}} \sum_{\alpha\geq \beta} C_{\alpha,\beta} \int_{\R^d} (\partial_\xi^{\alpha-\beta} q) (\partial_\xi^\beta e^{it\mu_{\mathcal{A}}} ) e^{i y\cdot \xi} \dd\xi\\
& =  \frac{i^{|\alpha|}}{(2\pi)^\frac{d}{2}} \sum_{\alpha\geq \beta} C_{\alpha,\beta}\big (  \int_{\R^d} \chi_R \cdots \dd \xi + \int_{\R^d} (1-\chi_R ) \cdots \dd\xi \big ).
\end{align*}
Here $\chi_R$ is he cut-off function as in the proof of \eqref{proof.lem.G_p.2}. By the definition of $q$ we see 
\begin{align*}
|\int_{\R^d} \chi_R \cdots \dd \xi | \leq C \int_{|\xi|\leq 2 R} |\xi |^{-|\alpha | +1} (1+|\xi|)^{-1-\epsilon} \dd\xi \leq C R^{d+1-|\alpha| },
\end{align*}
and 
\begin{align*}
|y  \int_{\R^d} (1-\chi_R ) \cdots \dd\xi |\leq C \int_{|\xi|\geq R} |\xi |^{-|\alpha |} (1+|\xi| )^{ - 1 -\epsilon}\dd\xi\leq C R^{d-|\alpha| -1-\epsilon}.
\end{align*} 
Thus it follows that $|G_q(x,y,t) |\leq C |y|^{-|\alpha|} (R^{d+1-|\alpha|} + |y|^{-1} R^{d-|\alpha|-1-\epsilon} )$ for $x,y,\in \R^d$ and $0<t<2$. If $|y|\leq 1$ then we take $|\alpha|=d-1$, while if $|y|>1$ then take $|\alpha|=d$. Then, putting  $R=|y|^{-\kappa}$ with sufficiently small $\kappa>0$, we get \eqref{proof.lem.G_p.3}. The proof is complete.

\end{document}